\documentclass[onefignum,onetabnum]{siamonline220329}

\usepackage{lipsum}
\usepackage{amsfonts}
\usepackage{graphicx}
\usepackage{epstopdf}
\usepackage{algorithmic}
\usepackage{multirow}
\usepackage{caption}

\usepackage{pgfplots}
\pgfplotsset{compat=newest}
\usepackage{tikz}
\usetikzlibrary{matrix} 

\ifpdf
  \DeclareGraphicsExtensions{.eps,.pdf,.png,.jpg}
\else
  \DeclareGraphicsExtensions{.eps}
\fi

\usepackage{xcolor}

\usepackage{enumitem}
\setlist[enumerate]{leftmargin=.5in}
\setlist[itemize]{leftmargin=.5in}

\newcommand{\pl}[1]{{#1}^{+}}
\newcommand{\T}{^{\hthinsp\mathrm{\scriptstyle T}}}
\newcommand{\Ti}{^{\hthinsp\mathrm{\scriptstyle -T}}}
\newcommand{\hthinsp}{\mskip  1   mu}    
\newcommand{\extsup}[3]{{#1}^{\textnormal{#2}}_{#3}}
\newcommand{\clb}[1]{{\color{blue} #1}}
\newcommand{\bmat}[1]{\begin{bmatrix}#1\end{bmatrix}} %

\newsiamremark{remark}{Remark}
\newsiamremark{hypothesis}{Hypothesis}
\crefname{hypothesis}{Hypothesis}{Hypotheses}
\newsiamthm{claim}{Claim}

\headers{Learning the Bidiagonalization}{J. J. Brust, and M. A. Saunders}

\title{Fast and Accurate SVD-Type Updating in Streaming Data\thanks{Submitted to the editors Fall 2025.  Version of \today.
\funding{This work was partially funded by the startup fund at Arizona State University under Grant PG16270.}}}

\author{Johannes J. Brust\thanks{School of Mathematical and Statistical Sciences, Arizona State University, Tempe, AZ 
  (\email{jjbrust@asu.edu}).}
  \and Michael A. Saunders\thanks{%
  Dept of Management Science and Engineering, Stanford University, Stanford, CA (\email{saunders@stanford.edu}).
  }}

\usepackage{amsopn}

\usepackage{comment,xspace}

\ifpdf
\hypersetup{
  pdftitle={SVD-Type Updates},
  pdfauthor={J. J. Brust, and M. A. Saunders}
}
\fi

\begin{document}

\maketitle

\begin{abstract}
For a datastream, the change over a short interval is often of low rank. For
high throughput information arranged in matrix format, recomputing an optimal
SVD approximation after each step is typically prohibitive. Instead, incremental and truncated updating strategies are used, which may not scale for large truncation ranks. Therefore, we 
propose a set of efficient new algorithms that update a bidiagonal factorization, and which are similarly accurate as the SVD methods. In particular, we develop a compact 
Householder-type algorithm that decouples a sparse part from a low-rank
update and has about half the memory requirements of standard bidiagonalization methods.
A second algorithm based on Givens rotations has only about 10 flops 
per rotation and scales quadratically with the problem size, compared
to a typical cubic scaling. The algorithm
is therefore effective for processing high-throughput updates, as 
we demonstrate in tracking large subspaces of recommendation systems and networks,
and when compared to well known software such as LAPACK or the incremental SVD.
\end{abstract}

\begin{keywords}
Streaming data, Low-rank updates, SVD, Bidiagonal decomposition, Recommendation system, Subspace tracking,
incremental SVD
\end{keywords}

\begin{MSCcodes}
 65F55, 68T09, 65Y20
\end{MSCcodes}

\section{Introduction}
Data compression, representation and updating are essential for 
applications in machine learning, data science, imaging
and many more (Sun et al.\ \cite{SGLTU20} and recent work \cite{jeong2025stochastic,cohen2025efficient,saibaba2025randomized,BrBr23}). We suppose the data is 
represented in a dense or sparse matrix $A \in \mathbb{R}^{m \times n}$. For lossy
compression of (for example) image or video data, the truncated Singular Value
Decomposition (SVD) produces the best rank-$r$ approximation of the 
data (Eckart and Young \cite{eckart1936approximation}). For streaming 
applications such as recommendation systems, a subspace tracking method
based on a truncated SVD can be 
effective (Brandt \cite{brand2006fast,brand2002incremental}). Typically only a small
``latent'' subspace for a large database of user choices is sufficient
to capture the essential preferences for many users (Sarwar et al.\ \cite{sarwar2002incremental}). 
Recently, Deng et al.~\cite{dengIncSVD} proposed a performant subspace tracking
method for representation learning.
For optimization, updates to Hessian approximations
are essential for algorithms to
incorporate the most recent information (Brust \cite{brust2022large}, 
Brust and Gill \cite{brust2023ldl}). A useful representation that captures
many applications is the low-rank update model
\begin{equation}
\label{eq:main}
    \pl{A} = A + B C\T,
\end{equation}
where $ \pl{A} \in \mathbb{R}^{m \times n}$ is the updated system,
$A \in \mathbb{R}^{m \times n} $ is the current data, and the matrices
$ B \in \mathbb{R}^{m \times r}  $ and $ C \in \mathbb{R}^{n \times r}  $ contain new information.
For instance, if $A$ is a matrix of user ratings (say, movie ratings), $B=-Ae$ ($e$ being the vector of all ones), and $C=e/n$, then 
$\pl{A}$ is the mean adjusted ratings data. Alternatively, for the same $A$,
if $C$ is a column of the identity matrix and $B$ is a column of preference
data for a new user, the resulting $\pl{A}$ represents an augmented matrix/table
that includes the new user. Update \cref{eq:main} uses $mnr$ multiplications. 

Often it is important not only to compute $\pl{A}$
but to use it for some analysis. Therefore, $\pl{A}$ is typically further factored
to obtain a useful representation of the data. 
Depending on the application, high accuracy may or may not be desired. However, 
it is vital that new
information be incorporated efficiently, that sparsity be maintained, and
that memory be kept low (especially for large problems). Downstream, the SVD
is used for many data analysis tasks and is itself a two-phase algorithm. 
First, the matrix is reduced to bidiagonal 
form by a sequence of finite matrix operations.  The bidiagonal factorization
of real $A$ always exists: 
\begin{equation}
\label{eq:bidiag}
A = Q B P\T,
\end{equation}
where $Q \in \mathbb{R}^{m \times m}$ and $P \in \mathbb{R}^{n \times n}$ are orthogonal
and $ B \in \mathbb{R}^{m \times n} $ is upper bidiagonal. 
Second,
the singular values and 
vectors are found by iteratively reducing
the bidiagonal to diagonal, with typically fast convergence for each singular
value (Watkins \cite{watkins2004fundamentals}).

\subsection{Contributions}
\label{sec:contrib}
We develop new methods to update the bidiagonal factorization when low-rank updates are added to a previous bidiagonal factorization. The methods do not need 
a refactorization, and they maintain the data sparsity. In particular, they
can be implemented as \emph{matrix-free} algorithms that only need access to matrix-vector
products, or previously computed factors. Because bidiagonalization is the main step in an          
SVD, the new methods can be used for data compression, data analysis, or data approximation.
We demonstrate how our results are effective for handling large high-throughput streaming data,
including when compared to widely used existing methods.

\subsection{Notation}
\label{sec:notation}
We use Householder notation: upper-case roman letters for matrices, lower-case roman letters for vectors, and
lower-case Greek letters for scalars. An exception is the diagonal matrix of 
singular values $ \Sigma = \textnormal{diag}(\sigma_i), \: \sigma_i \ge \sigma_{i+1} $
from the SVD $A = U \Sigma V\T $. The integer $k\ge 0$ denotes 
the iteration index. We label the minimum dimension by $ t = \textnormal{min}(m,n)$. A bidiagonal matrix has a superdiagonal band:
\begin{equation*}
B = 
\begin{bmatrix}
   \alpha_1 & \beta_1 & &
\\          & \alpha_2 & \ddots &
\\          &          & \ddots & \beta_{t-1}
\\[4pt]     &          &        & \alpha_t
\\[20pt]
\end{bmatrix}.
\end{equation*}
By definition, $\beta_0 = 0$. 
We denote the identity matrix by $I$ with dimensions that depend on the context. 
The Euclidean norm is abbreviated by $ \| \cdot \|_2 := \| \cdot \| $ and the Frobenius norm is denoted by $ \| \cdot \|_F$.
We index arrays and vectors using colon notation such as $B(:,k) = b_k $ and $b_k(k\!-\!1\!:\!k) = [\beta_{k-1}, \: \: \alpha_k]\T  $.

\subsection{Motivation}
\label{sec:motiv}
It is well established that the best rank-$r$ approximation of a matrix $A$ is given by the truncated
SVD:
\begin{equation}
\label{eq:trSVD}
A\approx \extsup{A}{SVD}{r} = \sum_{i=1}^r \sigma_i u_i v\T_i.
\end{equation}
The approximation error can be computed exactly (Eckart and Young \cite{eckart1936approximation}) as the sum of truncated squared singular values:
\begin{equation}
\label{eq:errtrSVD}
\| A - \extsup{A}{SVD}{r} \|^2_{F} = \sum_{i=r+1}^t \sigma^2_i.
\end{equation} 
Hence $r=t=\textnormal{min}(m,n)$
results in an exact reconstruction of the data. For practical applications, this ideal value of $r$ is typically 
too costly. Even forming the truncated SVD is costly, because a full factorization is needed
to find the $r$ largest singular values. 
An alternative to truncated SVD
is the truncated bidiagonal decomposition (BD).
Similar to \cref{eq:trSVD}, a rank-$r$ approximation of $A$ is obtained from \eqref{eq:bidiag} as
\begin{equation}
\label{eq:trBD}
A\approx \extsup{A}{BD}{r} 
= \sum_{i=1}^r \alpha_i q_i p\T_i + \beta_{i-1} q_{i-1} p\T_i
= \sum_{i=1}^r (\alpha_i q_i + \beta_{i-1} q_{i-1}) p\T_i.
\end{equation}
As before, the approximation error can be found exactly (see Supplement \Cref{sec:AppBidiagtrErr}): 
\begin{equation}
\label{eq:errtrBD}
\| A - \extsup{A}{BD}{r} \|^2_{F} = \sum_{i=r+1}^t \alpha^2_i + \beta^2_{i-1}.
\end{equation}
The best BD approximation can therefore be found by keeping
the $r$ pairs $\{\alpha_i,\beta_{i-1}\}$ for which $\alpha^2_{i} + \beta^2_{i-1}$ is largest. 
If the $\beta_i$'s are zero, the bidiagonal decomposition is the same as the SVD.
Lower and upper bounds for the difference between the SVD and bidiagonal approximations are
\begin{equation}
    \label{eq:diffErr}
    \sum_{i=1}^r \sigma^2_i - (\alpha^2_i + \beta^2_{i-1}) \le \| \extsup{A}{SVD}{r} - \extsup{A}{BD}{r} \|^2_F \le \sum_{i=i_j}^{t_j} \sigma^2_i + \alpha^2_i + \beta^2_{i-1}, 
\end{equation}
\vspace{-20pt}
\begin{equation*}
i_j \in \{1,r+1\}, \quad t_j \in \{r, t\}, \quad j \in \{1,2\}.
\end{equation*}
To understand the upper bound \cref{eq:diffErr} for $j=2$, namely $i_2=r+1$ and $t_2 = t = \textnormal{min}(m,n)$,
observe that the middle term can be written as $ \| (A-\extsup{A}{SVD}{r}) - (A-\extsup{A}{BD}{r}) \|^2_F  $. 
Applying the triangle inequality and \cref{eq:errtrSVD,eq:errtrBD} yields one upper bound. A derivation for $j=1$
and the lower bound is in \Cref{sec:AppBetweenErr} (Supplementary Materials). 

The bounds are illustrated in \cref{fig:bounds_supp} in the Supplementary Materials.

\begin{figure}   

\includegraphics[width=0.49\textwidth,trim={0 0 2.6cm 0cm},clip]{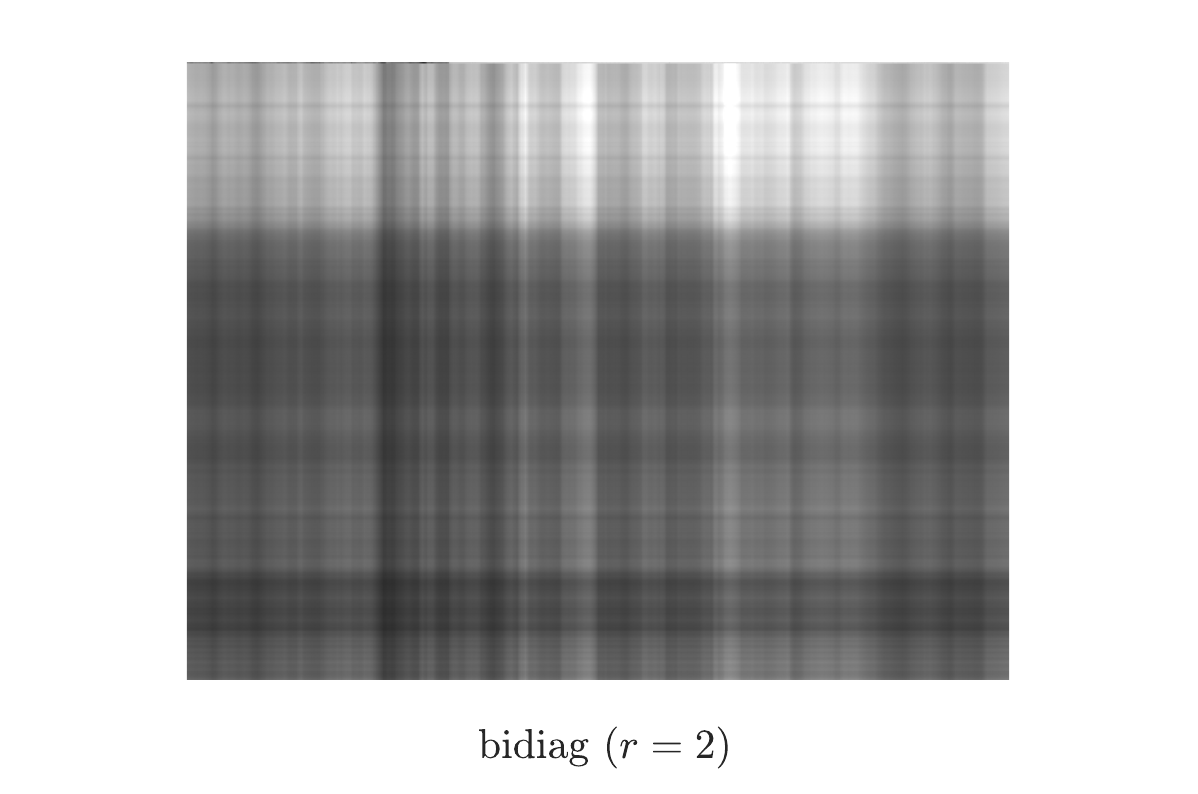}%
\includegraphics[width=0.49\textwidth,trim={2.6cm 0 0 1.4cm},clip]{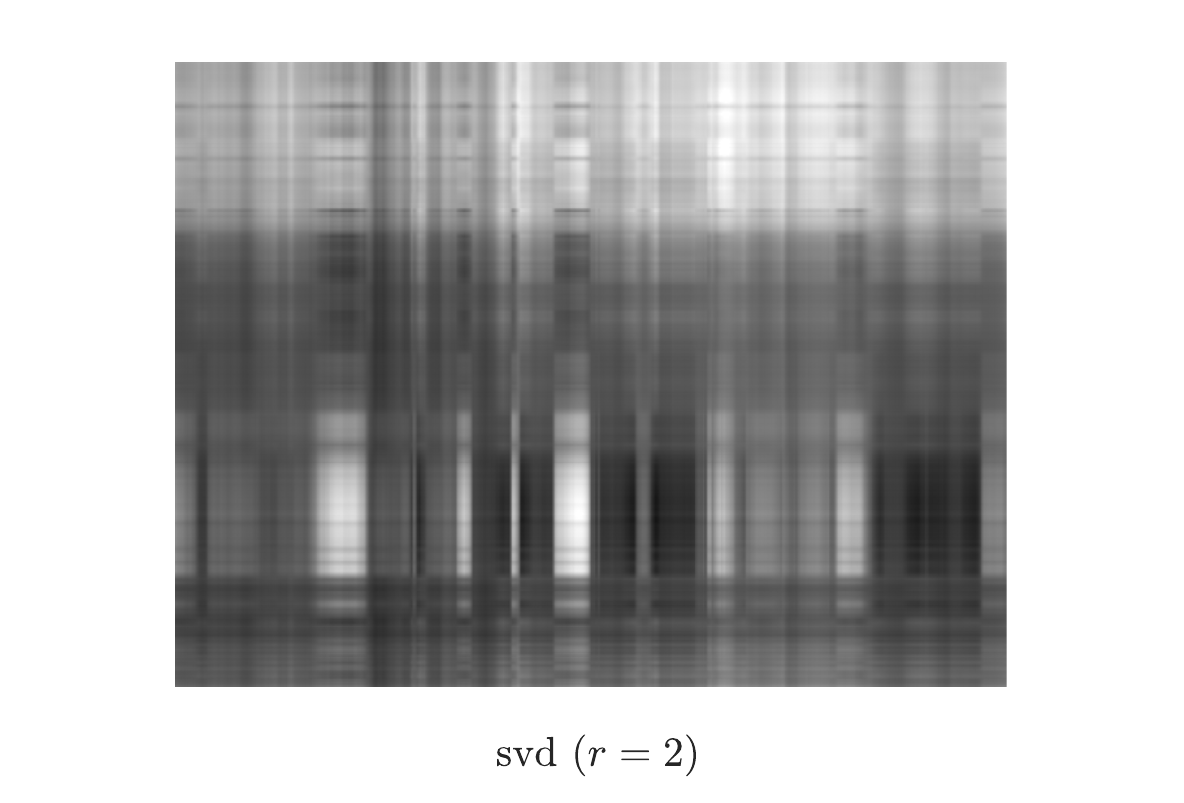}
\includegraphics[width=0.49\textwidth,trim={0 0 2.5cm 0cm},clip]{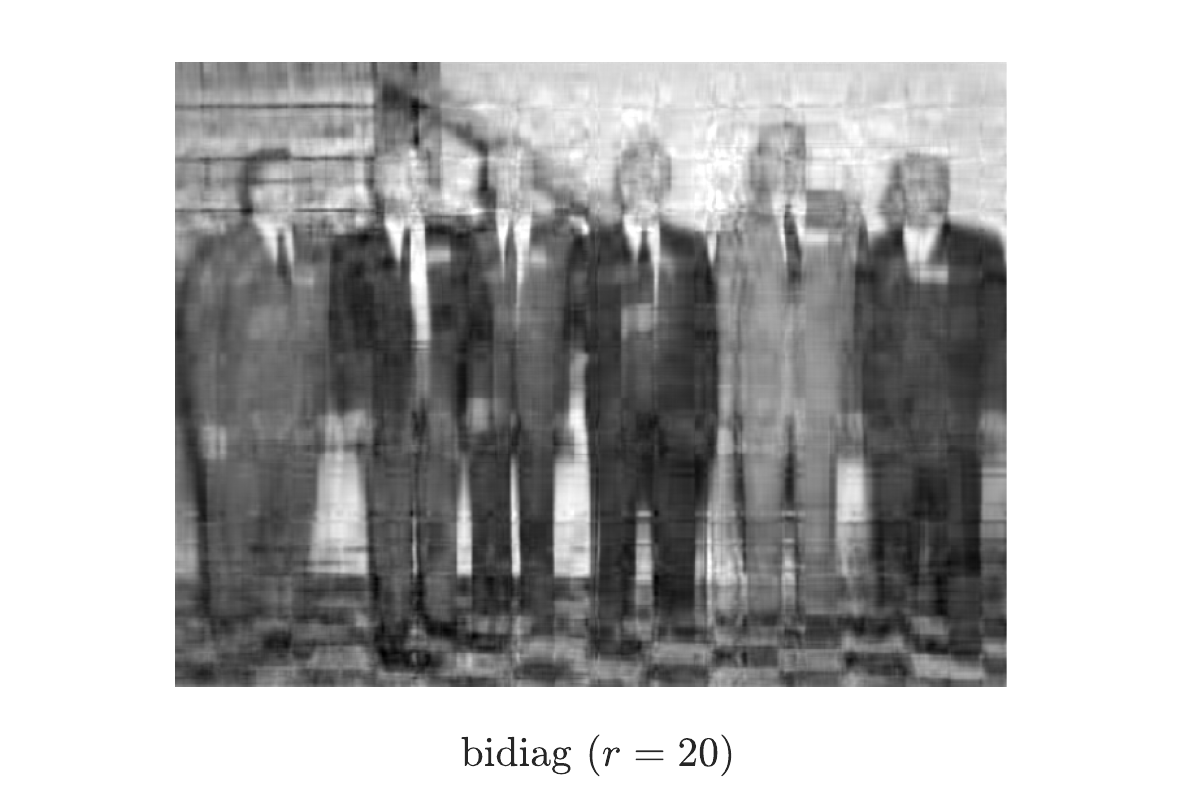}%
\includegraphics[width=0.49\textwidth,trim={2.5cm 0 0 0cm},clip]{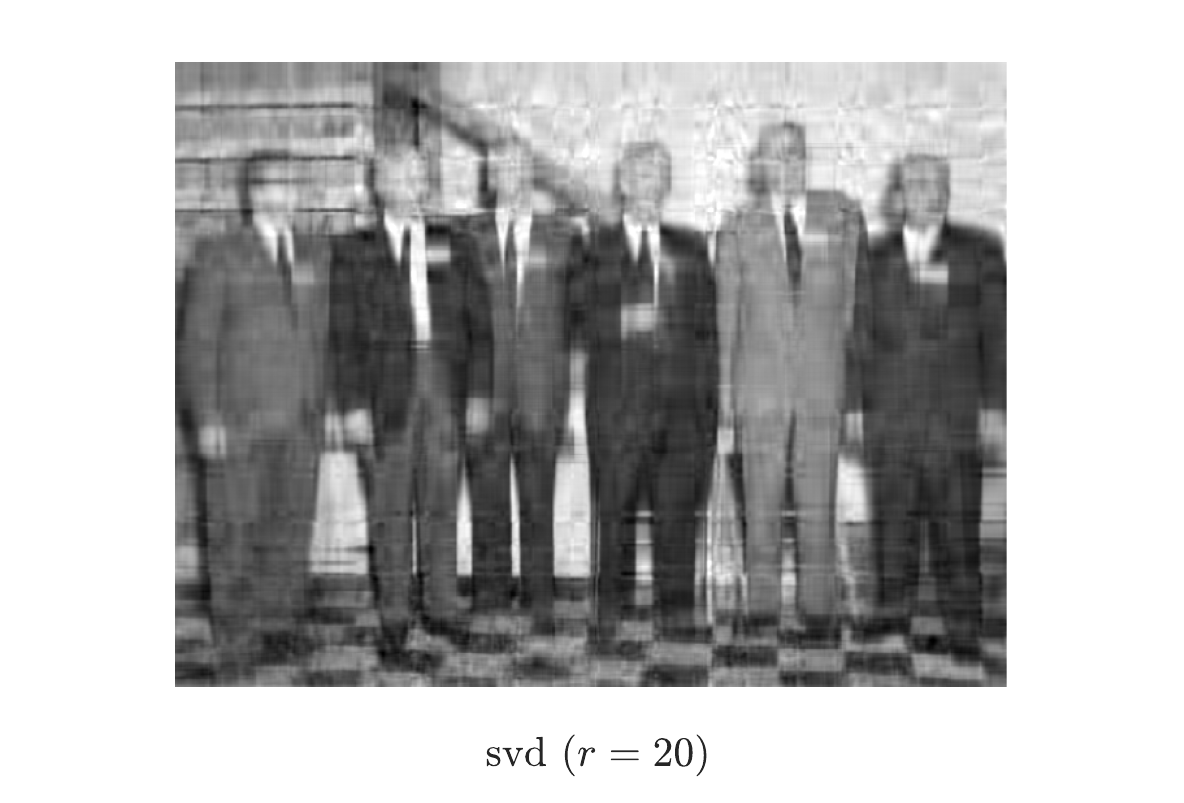}
\includegraphics[width=0.49\textwidth,trim={0 0 2.5cm 0cm},clip]{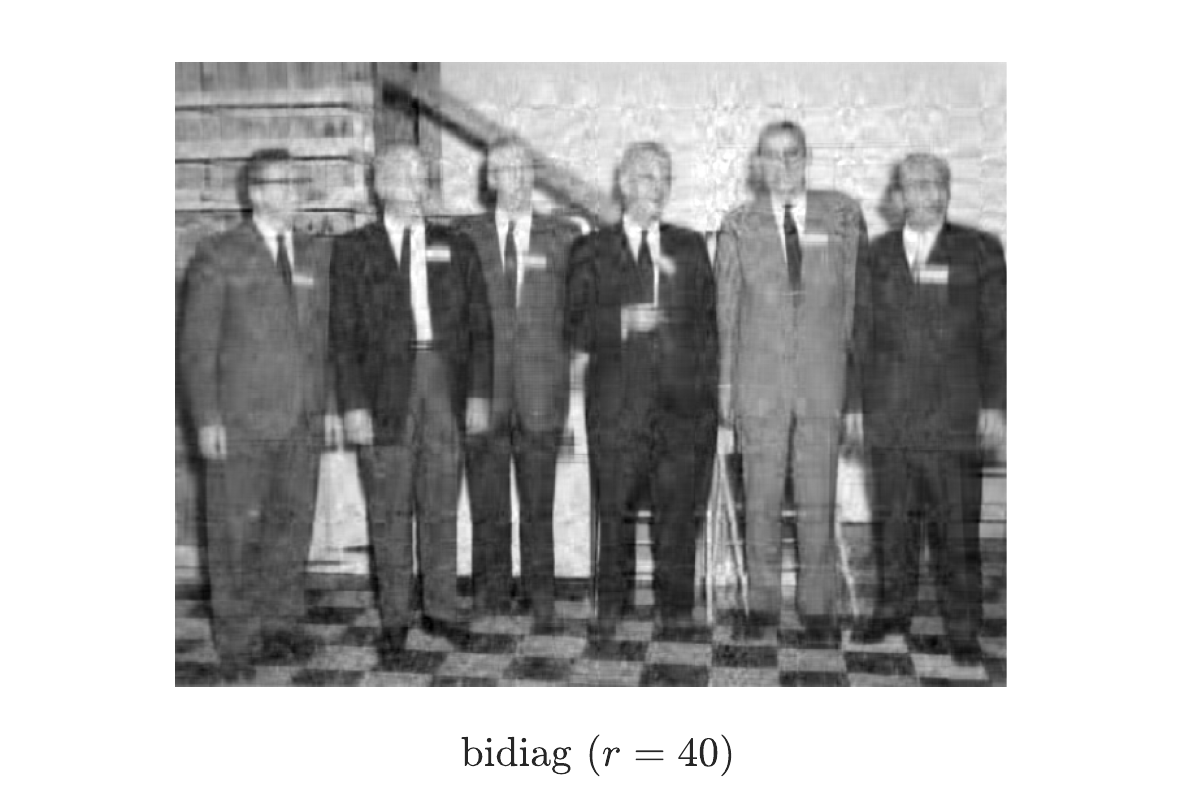}%
\includegraphics[width=0.49\textwidth,trim={2.5cm 0 0 0cm},clip]{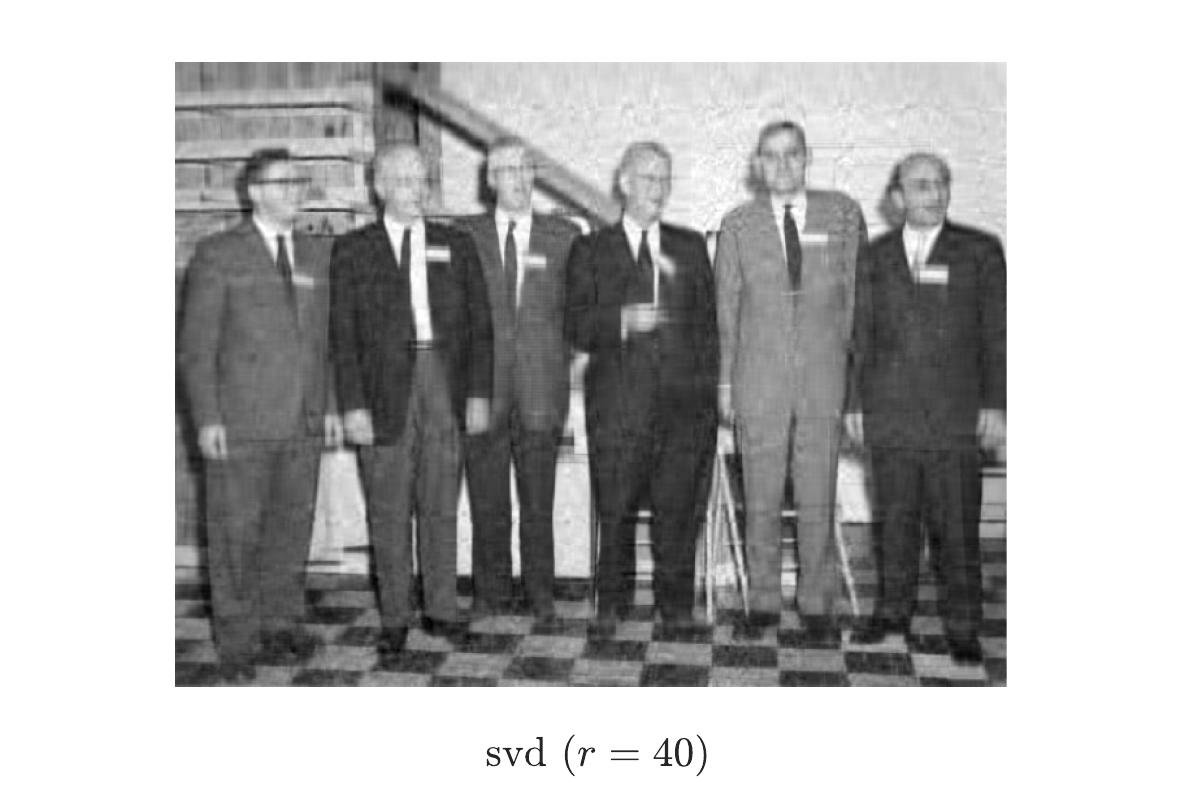}

\caption{Low-rank approximations of an image using truncated BD and truncated SVD. 
The BD is qualitatively similar to the SVD.}
\label{fig:reconstruct}
\end{figure}

For practical purposes the bidiagonal decomposition
can give useful reconstructions (Simon and Zha \cite{simon2000low}). As an example, \cref{fig:reconstruct} shows the 
reconstruction of a $359\times 371$ grayscale image using a sequence of low-rank approximations 
based on truncated BD or SVD. 

\section{Methods}
\label{sec:methods}

\subsection{Bidiagonalization}
\label{sec:bidiag}
The bidiagonal factorization can be computed in different ways. The most 
 stable approach uses alternating sequences of orthogonal transformations on the left and right, as implemented in LAPACK (Anderson et al.\ \cite{anderson1999lapack}). 
 It is the first step in the Golub-Reinsch algorithm for computing the SVD of a dense matrix \cite{Golub:1970:SVD}.

For sparse $A$,
the Golub-Kahan bidiagonalization (GKB)
\cite{GK65}
is an iterative method based on two short recurrences.

Although there seems to be no existing work on updating a bidiagonal 
factorization with low-rank updates, Brand \cite{brand2006fast,brand2002incremental} has proposed a method for updating a thin SVD. The approach is based on an SVD of a diagonal plus low-rank matrix, which is easier than the original problem. 
Other SVD updating algorithms are given by Moonen et al.\ \cite{moonen1992singular}, Zha and Simon~\cite{zha1999updating}, 
Vecharynski and Saad~\cite{vecharynski2014fast}, and Yamazaki et al.~\cite{yamazaki2015randomized}. Recently,
Deng et al.~\cite{dengIncSVD} compared methods \cite{zha1999updating,vecharynski2014fast,yamazaki2015randomized} and their own
for large machine learning settings.
In the following, we 
describe known and new techniques for computing a bidiagonal factorization.

\subsection{Householder Reflectors}
\label{sec:house}
The most stable method to implement bidiagonalization is to apply Householder reflectors
\cite{golub2013matrix,householder1958unitary} alternately on the left and right of $A$. Each reflector is an orthogonal symmetric rank-1 update to the identity:
\begin{equation}
    \label{eq:house}
    H_k = I - \tau_k y_k y_k\T, \quad \tau_k = \frac{2}{ \| y_k \|^2_2 }.
\end{equation}
The vector $y_k$ is the ``essential part'' of the reflector and is chosen
to create zeros in a column or row of $A$. Two choices of $y_k$ that
zero all elements of a vector $ a_k $ below $a_{k}(k)$
in the product $H_k a_k$ are $ y^{\textnormal{left}}_k := y^{\textnormal{l}}_k =
\begin{bmatrix} 0(\!1\!:\!k\!-\!1) \\ a_k(\!k\!:\!m) \end{bmatrix} \pm \| a_k(\!k\!:\!m) \| e_k $. Many software implementations additionally scale $y_k$  so that 
$y_k(k)=1$ to economize on storage. They may also 
ensure that 
$ \| y_k \|_2 = 1 $. If $a^*_k$ denotes the $k^{\textnormal{th}}$ row of $A$,
a Householder reflector can be used to zero out all elements to
the right of element $k+1$. The reflectors in this case are defined 
by $ (y^{\textnormal{right}}_{k+1})\T := (y^{\textnormal{r}}_{k+1})\T =
\begin{bmatrix} 0(\!1\!:\!k\!) & a^*_k(\!k\!+\!1\!:\!n) \end{bmatrix} \pm \| a^*_k(\!k\!+1\!:\!n) \| e_{k+1}\T $.
For illustration, after the $k^{\textnormal{th}}$ step of a bidiagonalization 
with Householder reflectors, the nonzero pattern looks like
\begin{equation}
\label{eq:housesparse}
\extsup{H}{l}{k}\extsup{H}{l}{k-1}\dots \extsup{H}{l}{1} A \extsup{H}{r}{2} \dots \extsup{H}{r}{k+1}
= 
\left[
\begin{array}{c c c c | c c c}
    \alpha_1    & \beta_1   &           & & & & \\
                & \ddots    & \ddots        &               & & & \\
                &           & \alpha_{k-1}  & \beta_{k-1}   &        &          & \\
\hline                &           &               & \alpha_k      & \times & \cdots   & \times \\
                &           &               &               & \times & \ddots   & \times \\
                &           &               &               & \vdots & \ddots   & \vdots \\
                &           &               &               & \times & \cdots   & \times \\
\end{array}
\right],
\end{equation}
where $ \extsup{H}{l}{i} $ and $ \extsup{H}{r}{i+1} $ are reflectors for zeroing columns and rows respectively,
and ``$\times$'' represents an arbitrary element. When the process ends, the factorization is upper bidiagonal:
\begin{equation*}
    \extsup{H}{l}{t-1}\cdots \extsup{H}{l}{k} \cdots \extsup{H}{l}{1} A \extsup{H}{r}{2} \cdots \extsup{H}{r}{k} \cdots \extsup{H}{r}{t-1} 
    := Q\T A P = B.
\end{equation*}
Householder bidiagonalization is well suited to dense matrices, but 
does not exploit sparse or structured $A$. 
At iteration $k$, a reflector  $\extsup{H}{l}{k}a_k$ applied from the left costs
$(m-k+1)(n-k+1)$ multiplications, and the product $a^*_k \extsup{H}{r}{k+1}$ from the right has $(m-k)(n-k)$ multiplications.
For $t=n\le m$, the total multiplication count is
\begin{equation}
\label{eq:cplxhouse}
\sum_{k=1}^{t-1}(m-k+1)(n-k+1) + (m-k)(n-k) \le \sum_{k=1}^{n-1}2(m-k+1)(n-k+1) =  \mathcal{O} \big( (mn-\frac{1}{3}n^2) t \big). %
\end{equation}
Thus, the total factorization cost scales as $mnt$ for large $t$. The same analysis with $ t = m \le n$ gives a symmetric result of $ \mathcal{O} \big( (nm-\frac{1}{3}m^2) t \big ) $.

\subsection{Golub-Kahan Bidiagonalization}
\label{sec:gkb}

GKB \cite{GK65} is based on the existence of the factorization 
$$ 
Q B P\T = A 
$$ 
and
the corresponding equations
\begin{align}
\label{eq:qbap} QB &= AP \\
\label{eq:pbaq} PB\T &= A\T Q.
\end{align}
Multiplying \cref{eq:qbap,eq:pbaq} by $e_k$ from the right on both sides yields the recurrence relations
\begin{align}
   \label{eq:recp} \beta_{k-1} q_{k-1} + \alpha_k q_k &= Ap_k,
\\ \label{eq:recq} \alpha_k p_k + \beta_k p_{k+1} &= A\T q_k.
\end{align}
After \cref{eq:recp}--\cref{eq:recq} are initialized by an arbitrary unit-length vector $p_1$,
all subsequent $\alpha_k, \beta_k, q_k$ and $p_k$ can be found from the relations.
Note that $A$ is not overwritten, as it is in Householder bidiagonalization. In exact arithmetic the sparsity pattern at iteration $k$ is
\begin{equation*}
    Q_k \T A P_k = 
    \begin{bmatrix}
        \alpha_1    & \beta_ 1  &   & \\
                    & \ddots    & \ddots & \\
                    &           & \alpha_{k-1}  & \beta_{k-1} \\
                    &           &               & \alpha_k
    \end{bmatrix},
\end{equation*}
and the method does not carry additional nonzero elements like in \cref{eq:housesparse}.
The computational complexity mainly depends on the cost of matrix-vector products.
Typically, $A$ is sparse for GKB and the products are inexpensive. 
To be general, we assume for the moment that a product with $A$ or $A\T\!$ costs $mn$ multiplications.
Since there are two such multiplications in \cref{eq:recp,eq:recq}, the dominant 
operations are of order $2 m n t$ multiplications.

GKB is suited to large sparse systems because it only needs matrix-vector products and does
not overwrite the original system. 
A drawback is the loss of orthogonality,
which diminishes the method's accuracy. A variety of reorthgonalization strategies exist, such as Simon and Zha \cite{simon2000low}. For rectangular data with
one dimension much larger than the other, Simon and Zha orthogonalize only the
low-dimensional vectors (short space) because the orthogonality of the large vectors (long space)
is related to the errors in the short space. This strategy can be computationally effective
for very nonsquare data.




\section{New Methods}
\label{sec:newmethods}
\Cref{sec:methods} describes existing methods. We now introduce
three new methods.

\subsection{Randomized Bidiagonal Decomposition}
\label{sec:rBD}
We first propose a randomized bidiagonal factorization similar to the 
randomized SVD (RSVD) \cite{musco2015randomized}, which we call 
randomized bidiagonal decomposition (RBD). RBD has advantages
similar to RSVD, namely that the approximation error $ \| A - \extsup{A}{RBD}{r} \|^2_F $ 
is close to the best error $ \sum_{i=r+1}^t \sigma^2_i $. As for RSVD, RBD may perform suboptimally when the $\alpha_i$ and $\beta_i$'s do not decay
quickly to zero. A desirable feature of RBD is that it can be succinctly
summarized. Let $ S $ be a random sketching matrix, and consider this 
sequence of assignments:

\medskip

RBD (Randomized Bidiagonal Decomposition)

\medskip

\begin{tabular}{l l}
   $Y := AS$ &   \clb{// columnspace sketch using a random $S \in \mathbb{R}^{n \times r} $ }
\\ $ \extsup{Q}{Y}{r} R  := Y $ & \clb{// thin QR decomposition}
\\ $Z := (\extsup{Q}{Y}{r})\T A$ & \clb{// rowspace sketch}
\\ $\extsup{Q}{Z}{r} B_r P_r\T := Z $ & \clb{// bidiagonal decomposition}
\\ $Q_r := \extsup{Q}{Y}{r}\extsup{Q}{Z}{r} $ & \clb{// orthonormal matrix}
\\ $ A_r := Q_r B_r P_r\T $ & \clb{// final approximation}
\end{tabular}


\bigskip

\noindent We may further implement the SVD of $B_r $ in a last step
$ \extsup{U}{B}{r} \extsup{\Sigma}{B}{r} {\extsup{P}{B}{r}} \T := B_r $, which
reduces RBD to the RSVD. However, the RBD gives the option of preprocessing the bidiagonal beforehand; for example,
by truncating small columns in $B_r$.

\subsection{Learning Low-Rank Updates}
\label{sec:lowrankupdates}
In an online setting, data is continuously updated. We consider rank-1 modifications,
which may be readily extended to higher ranks by repeating the process:
\begin{equation}
    \label{eq:rk1}
    \pl{A} = A + b c\T.
\end{equation}
We seek a new decomposition
\begin{equation}
\label{eq:bidiagp}
\pl{Q} \pl{B} {\pl{P}} \T = Q B P\T + b c\T
\end{equation}
by updating the previous decomposition. 
Orthogonality of $Q$ and $P$ gives
\begin{equation}
\label{eq:bidiagsubproblem}
    QBP\T + b c\T = Q ( B + Q\T b c\T P ) P\T := Q ( B + \pl{b} {\pl{c}}\T ) P\T,
\end{equation}
where $ \pl{b} := Q\T b $ and $ \pl{c} = P\T c $. Observe that 
$ B + \pl{b} {\pl{c}}\T $ is bidiagonal with a rank-1 modification, which should be a simpler problem than the original. If we can
decompose
\begin{equation*}
    B + \pl{b} {\pl{c}}\T = Q_1 \pl{B} P_1 \T,
\end{equation*}
the streaming update has the form
\begin{equation}
\label{eq:bidiagplus}
A + b c\T = Q( B + \pl{b} {\pl{c}}\T ) P\T =
            Q Q_1 \pl{B} P_1 \T P\T := 
            \pl{Q} \pl{B} {\pl{P}}\T.
\end{equation}
The last equality in \cref{eq:bidiagplus} is the desired new decomposition.
Observe that $ B + \pl{b} {\pl{c}} \T $ has a special sparsity pattern:
\begin{equation}
    \label{eq:sparsity}
    B + \pl{b} {\pl{c}} \T = 
    \left[
    \begin{array}{c c c c }
    \times    & \times   &           & \\
                & \ddots    & \ddots    & \\
                &           &           \times  & \times \\
                &           &           &               \times
    \end{array}
    \right] +
    \left[
    \begin{array}{c}
        \times \\
        \times \\
        \vdots \\
        \times
    \end{array}
    \right]
    \left[
        \begin{array}{c c c c}
            \times & \times & \hdots & \times 
        \end{array}
    \right].
\end{equation}
Because decomposing \eqref{eq:sparsity} is an essential step of the method,
we consider a few alternatives next.

\subsection{Householder Low-Rank}
\label{sec:houselr}
An immediate approach is to apply Householder reflectors to the special bidiagonal update \cref{eq:sparsity}.
Applying the first reflector from the left zeroes all but the first element in column 1 (as desired).
However, it also creates fill-in within \cref{eq:sparsity}:
\begin{equation}
    \label{eq:bfill}
   \extsup{H}{l}{1} ( B + \pl{b} {\pl{c}} \T ) =
   \left[
    \begin{array}{c c c c }
    \times    & \times              &  \hdots           & \times \\
                & \vdots            & \vdots            & \vdots \\
                & \vdots            &           \vdots  & \vdots \\
                & \times            &         \hdots          &               \times
    \end{array}
    \right].
\end{equation}
For large matrices, $mn$ storage of nonzero elements is not feasible. We develop a Householder method
that does not create fill-in in the usual sense and keeps the zeros in the bidiagonal intact, by storing the rank-1 updates
incrementally.

For the remainder of this section we simplify the notation by denoting left Householder vectors by 
$ y_k $ and right ones by $ w_k $: 
\begin{equation}
    \label{eq:yw}
    y_k := y^{\textnormal{l}}_k \quad \text{ and } \quad  w_k :=  y^{\textnormal{r}}_k . 
\end{equation}
We assume the vectors are normalized ($ \| y_k \|_2 = 1 $, $ \| w_k \|_2 = 1 $)
and define matrices that collect all previous $ \{y_i,w_{i+1} \}_{i=1}^k $:
\begin{equation}
    \label{eq:YW}
    Y_k = \begin{bmatrix}
        y_1 & y_2 & \ldots & y_k
    \end{bmatrix} \in \mathbb{R}^{m \times k}, \quad \quad 
    W_k = \begin{bmatrix}
        w_2 & w_3 & \ldots & w_{k+1}
    \end{bmatrix} \in \mathbb{R}^{n \times k}.
\end{equation}
The columns of $Y_k$ and $W_k$ form two special upper triangular systems
\begin{equation}
    \label{eq:TR}
    \left( T_k \right)_{ij} = 
    \begin{cases}
        2 y_i\T y_j     &   i < j \\
        1               &   i = j \\
        0               & \textnormal{otherwise},
    \end{cases}
    \quad \quad 
    \left( R_k \right)_{ij} = 
    \begin{cases}
        2 w_{i+1}\T w_{j+1}     &   i < j \\
        1                       &   i = j \\
        0                       & \textnormal{otherwise}.
    \end{cases}
\end{equation}

Remarkably, Walker \cite{walker1988implementation} and later Puglisi \cite{puglisi1992modification} have
shown that the so-called compact WY representation of the product of Householder matrices described by 
Bischof and Van Loan \cite{bischof1987wy} can be given in an analytic formula. In particular,
for a sequence of reflectors \eqref{eq:house} (with normalized vectors), their product is an orthogonal
matrix that may be expressed in terms of \eqref{eq:YW} and \eqref{eq:TR}:
\begin{equation}
    \label{eq:QP}
    Q_k \equiv \prod_{i=1}^k H^{\textnormal{l}}_i = I - 2 Y_k T_k^{-1} Y_k\T, \qquad
    P_k \equiv \prod_{i=2}^{k+1} H^{\textnormal{r}}_i = I - 2 W_k R_k^{-1} W_k\T.
\end{equation}
Among other things this enables block algorithms that exploit matrix-matrix operations. The expressions 
in \eqref{eq:QP} are more storage-efficient than the original WY representation \cite{bischof1987wy} 
(which stores two large matrices, one of which is dense),
and explicitly define the elements in the middle matrix, in contrast to even an improved WY representation
\cite{schreiber1989storage}. The LAPACK  bidiagonalization algorithms and implementations (dgebrd and zgebrd)
also use a form of compact Householder reflectors 
closely related to the WY representation (see e.g., Dongarra et al.\ \cite{dongarra1989block}).
These implementations
compute and store multiple dense intermediate matrices.
However, by exploiting Level-3 BLAS operations, zgebrd is one of the most effective available algorithms. To maintain the bidiagonal form in \eqref{eq:sparsity} and avoid the fill-in from using the Householder reflectors, 
we develop a new compact representation for the product of compact reflectors \emph{and} a bidiagonal plus rank-1. 
(An approach for keeping the Householder updates separate from the initial matrix is also mentioned in Howell et al.\ 
\cite{howell2008cache}.)
The transformation of some $A$ 
\eqref{eq:housesparse} looks like
\begin{equation}
    \label{eq:arec}
    A_{kk} = \extsup{H}{l}{k} A_{k-1k-1} \extsup{H}{r}{k+1} = \extsup{H}{l}{k} \cdots \extsup{H}{l}{1} A \extsup{H}{r}{2} \cdots \extsup{H}{r}{k+1} = Q_k \T A P_k, 
\end{equation}
where $A_{kk}$ represents the initial matrix transformed by $k$ left and right orthogonal matrices. Below, we use the compact representations of $Q_k$ and $P_k$ from \eqref{eq:QP}. Note 
that if we 
immediately apply recurrence \eqref{eq:arec} to $ A = B + \pl{b} \smash{{\pl{c}}} \T $, the nonzero pattern in \eqref{eq:bfill}
shows that the product is dense and that any bidiagonal structure is lost. However, by 
exploiting the compact forms of $Q_k$ and $P_k$ we derive an analytic representation that preserves the sparsity, i.e., 
the bidiagonal part. 

\begin{theorem}
\label{thm:comp}
Applying $k$ Householder reflectors of the form \eqref{eq:house} from the left and right, with essential vectors stored in 
$ Y_k $ and $ W_k $ \eqref{eq:YW}, \eqref{eq:yw} 
to $A = B + b c\T $ gives the explicit 
decomposition 
\begin{equation}
    \label{eq:comp}
    A_{kk} = Q_k\T A P_k =  B - 
    \begin{bmatrix}
        b & Y_k & BW_k
    \end{bmatrix}
    \begin{bmatrix}
        -1              & 0_{1 \times k }               & c\T W_k \\
        Y_k\T b    & \frac{1}{2} T_k\T & Y_k\T B W_k \\
        0_{k \times 1}               &   0_{k \times k}  & \frac{1}{2} R_k
    \end{bmatrix}^{-1}
    \begin{bmatrix}
        c \T \\
        Y_k \T B \\
        W_k \T
    \end{bmatrix},
\end{equation}
where $ T_k $ and $ R_k $ are upper triangular \eqref{eq:TR}.
\end{theorem}
\begin{proof}
After $k$ orthogonal transformations 
we have
\begin{equation*}
    A_{kk} = Q_k \T A P_k = (I - 2 Y_k T^{-1}_k Y_k \T ) \T (B + b c\T ) (I - 2W_k R_k^{-1} W_k \T ).
\end{equation*}
Multiplying and factoring terms gives 
the equivalent 
expression 
\begin{equation*}
    A_{kk} = B - 2 
    \begin{bmatrix}
        b & Y_k & BW_k
    \end{bmatrix}
    \begin{bmatrix}
        -\frac{1}{2}              & 0_{1 \times k }               & c\T W_k R_k^{-1} \\
        T_k\Ti Y_k\T b    & T_k\Ti & -2 T_k \Ti  Y_k\T (B + b c\T ) W_k R_k^{-1} \\
        0_{k \times 1}               &   0_{k \times k}  & R_k^{-1}
    \end{bmatrix}
    \begin{bmatrix}
        c \T \\
        Y_k \T B \\
        W_k \T
    \end{bmatrix}.
\end{equation*}
We now apply a sequence of two block inverses, the first for the top $ k+1 \times k+1  $ block and 
the second for the full $ (1 + 2k) \times (1 + 2k) $ system:
\begin{align*}
A_{kk} &=  B - 2 
    \begin{bmatrix}
        b & Y_k & BW_k
    \end{bmatrix}
    \begin{bmatrix}
       \begin{bmatrix}
            -2 & 0_{1 \times k} \\
            2Y_k\T b & T_k \T
       \end{bmatrix}^{-1}                &  \begin{array}{c} 
                                                c\T W_k R^{-1}_k \\
                                                -2 T_k\Ti Y_k\T (B + b c\T) W_k R^{-1}_k
                                            \end{array}       \\
         \begin{array}{c c}
            0_{k \times 1} & 0_{k \times k}
         \end{array}                          & R^{-1}_k
    \end{bmatrix}
    \begin{bmatrix}
        c \T \\
        Y_k \T B \\
        W_k \T
    \end{bmatrix} \\
    &= 
    B - 2 
    \begin{bmatrix}
        b & Y_k & BW_k
    \end{bmatrix}
    \begin{bmatrix}
        -2              & 0_{1 \times k }               & 2c\T W_k \\
        2Y_k\T b    & T_k\T & 2Y_k\T B W_k \\
        0_{k \times 1}               &   0_{k \times k}  & R_k
    \end{bmatrix}^{-1}
    \begin{bmatrix}
        c \T \\
        Y_k \T B \\
        W_k \T
    \end{bmatrix}.
\end{align*}
Finally, factorizing a 2 from the middle matrix yields the desired expression 
\eqref{eq:comp}.
\end{proof}

\begin{figure}   
\begin{center}
\begin{tikzpicture}[scale=.4] %
  \def\a{8.5}
  \def\b{0.5}
  \def\c{\a/2-\b/2}
  \def\ch{\a/4-\b/4}
  \def\d{\c+\b}
  \def\ofst{0.4}
  \draw (0,0) rectangle (\a,\a);

  \draw (\b,0) -- (\b,\a); 
  \draw (\d,0) -- (\d,\a);

  \draw (0,\c) -- (\a,\c); 
  \draw (0,\c+\c) -- (\a,\c+\c);

\fill[fill=yellow, fill opacity=0.3] (\a,0) -- (\a,\c) -- (\d,\c) -- cycle;
\fill[fill=gray, fill opacity=0.3] (0,\a-\b) -- (0,\c) -- (\b,\c) -- (\b,\a-\b) -- cycle;
\fill[fill=violet, fill opacity=0.3] (0,\a) -- (0,\a-\b) -- (\b,\c+\c) -- (\b,\a) -- cycle;
\fill[fill=blue, fill opacity=0.3] (\b,\a-\b) -- (\b,\c) -- (\d,\c) -- cycle;
\fill[fill=purple, fill opacity=0.3] (\d,\c) -- (\a,\c) -- (\a,\a-\b) -- (\d,\a-\b) -- cycle;
\fill[fill=red, fill opacity=0.3] (\d,\a) -- (\a,\a) -- (\a,\a-\b) -- (\d,\a-\b) -- cycle;

\node at (\b/2,\a+\ofst) {$1$};
\node at (\b+\ch,\a+\ofst) {$k$};
\node at (\d+\ch,\a+\ofst) {$k$};

\node at (0-\ofst,\c+\c+\b/2) {$1$};
\node at (0-\ofst,\c+\ch) {$k$};
\node at (0-\ofst,\ch) {$k$};


\def\ras{\a+2}
\def\rae{\a+\a}

\draw (\ras,0) rectangle (\a+\c+2,\rae);
\draw (\ras,\rae) -- (\ras+\c,\rae-\c);
\draw (\ras,\b) -- (\ras+\c,\b);
\fill[fill=yellow, fill opacity=0.3] (\ras,\rae) -- (\ras+\c,\rae-\c) -- (\ras+\c,\rae) -- cycle;
\fill[fill=green, fill opacity=0.3] (\ras,\rae) -- (\ras,\b) -- (\ras +\c,\b) -- (\ras +\c,\rae-\c) -- cycle;
\fill[fill=gray, fill opacity=0.3]  (\ras,\b) -- (\ras+\c,\b) -- (\ras+\c,0) -- (\ras,0) -- cycle;

\node at (\ras-\ofst, \a + \ofst) {$m$};
\node at (\ras+\ch, \rae + \ofst) {$k$};
\node at (\ras+\ch, \c+ \ch) {$Y_k$};

\node at (\ras+\ch + \ofst + \ofst, \rae -\ch + \ofst + \ofst) {$\underset{\textnormal{w/o diagonal}}{R_k}$};

\def\rbs{\a+2+\c+2}
\def\rbe{\a+\c}

\draw (\rbs,0) rectangle (\a+2+\c+2+\c,\rbe);
\draw (\rbs,\rbe) -- (\rbs+\c,\rbe-\c);
\draw (\rbs,\b) -- (\rbs+\c,\b);

\fill[fill=blue, fill opacity=0.3] (\rbs,\rbe) -- (\rbs+\c,\rbe-\c) -- (\rbs+\c,\rbe) -- cycle;
\fill[fill=cyan, fill opacity=0.3] (\rbs,\rbe) -- (\rbs,\b) -- (\rbs +\c,\b) -- (\rbs +\c,\rbe-\c) -- cycle;
\fill[fill=red, fill opacity=0.3]  (\rbs,\b) -- (\rbs+\c,\b) -- (\rbs+\c,0) -- (\rbs,0) -- cycle;

\node at (\rbs-\ofst, \c + \ch + \ofst) {$n$};
\node at (\rbs+\ch, \rbe + \ofst) {$k$};
\node at (\rbs+\ch, \c+ \ch) {$W_k$};

\node at (\rbs+\ch + \ofst + \ofst, \rbe - \ch + \ofst + \ofst) {$T_k$};

  
\end{tikzpicture}
\end{center}
\caption{Pattern of the middle matrix in Theorem 1, and the memory layout of the representation. Solving with the middle matrix is equivalent to a triangular solve. Importantly, the middle matrix
is never explicitly formed because all necessary information can be stored in two arrays of size $ m\times k  $ and $ n \times k $ and two $k$-vectors.}
\label{fig:mid}
\end{figure}
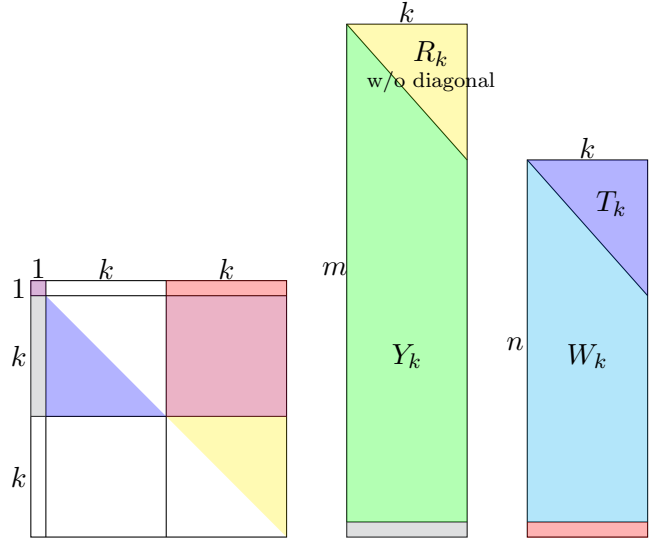

We can now use \cref{thm:comp} to deduce the bidiagonalization of $ B + \pl{b} {\pl{c}}\T $. Specifically,
after $n-1$ left and $n-2$ right reflectors we obtain the new bidiagonal $ \pl{B} = A_{n-1,n-2}  $.
It is enough to emphasize some properties of the representation in \eqref{eq:comp}. First, it is clear
that \cref{thm:comp} decouples the partial bidiagonal $A_{kk}$ into the sum of the previous 
bidiagonal $B$ and a rank-$(2k+1)$ correction. Therefore, there is an obvious separation into a sparse
and a dense part. Second, $A_{kk}$ can be implemented in a matrix-free way: 
$ A_{kk} $ is not explicitly formed but instead kept in factored form. The components
that define this matrix are $ B, W_k, Y_k, \pl{b} $ and $\pl{c}$. Further, the storage 
required for \eqref{eq:comp} (not counting $B$, $\pl{b}$ and $\pl{c}$ at $3n + m -1$ locations) is
$ mk + nk + 2k = (m+n+2)k $. We store the upper trapezoidal $Y_k$ and $W_k$ and the triangular 
$T_k$ and $R_k$ and two $k$-vectors. In fact, the essential nonzeros of $T_k$ fit into the zeros of $Y_k$, and so does
$R_k$ fit into $W_k$. We typically also store the two $k$-vectors $ Y_k \T \pl{b} $ and 
$ {\pl{c}}\T W_k $. In contrast, a standard bidiagonalization algorithm, like dgebrd \cite{anderson1999lapack}, requires
storage of a full dense $A_{kk}$, two upper trapezoidal and two dense $m\times k$
and $n \times k$ matrices. The total is $ (2m + 2n - k +1)k + mn $, which is more than twice the 
amount for \eqref{eq:comp}. Moreover, a standard bidiagonalization creates fill-in so that even
after the first iteration, $m \times n$ memory is needed. In contrast, 
the representation in \cref{thm:comp} has memory requirements that grow with iterations,
and stopping the factorization early results in memory used only up to that point. Such an 
advantage has also been noted in \cite{howell2008cache}.
Finally, solving with the middle matrix in \eqref{eq:comp} is efficient because it is equivalent to solving with a $ (1+2k) \times (1+2k) $ triangular matrix. \cref{fig:mid} illustrates 
the sparsity pattern and the memory layout of the representation.

To describe a practical algorithm we denote representation \eqref{eq:comp} by
\begin{equation}
    \label{eq:compshort}
    A_{kk} = Q_k \T (B + \pl{b} {\pl{c}}\T ) P_k = B - U_{kk} M^{-1}_{kk} V^T_{kk},
\end{equation}
where $ M_{kk} $, $ U_{kk} $ and $ V_{kk} $ are the middle, left and right matrices in the 
representation. \\

Algorithm 1 (BHU: Bidiagonal Householder Update) 

$ U_{00} = \pl{b}  $, $M_{00} = 1, V_{00} = \pl{c}, Y_{0} = [~], W_{0} = [~]  $ \\ 

\begin{tabular}{c c c}
\multicolumn{3}{l}{$\textnormal{for } k=1,\ldots,n-1 $} \\
    & \multicolumn{2}{l}{$ a_k = (B - U_{k-1k-1} M^{-1}_{k-1k-1} V_{k-1k-1} \T ) e_k  $} \\
    & \multicolumn{2}{l}{$ \alpha_k, y_k = \textnormal{house}(a_k)  $} \\
    & \multicolumn{2}{l}{$ Y_k = \begin{bmatrix} Y_{k-1} & y_k \end{bmatrix}  $} \\
    & \multicolumn{2}{l}{$ \textnormal{if } k = n-1 \textnormal{ then } \textnormal{ exit } $} \\
    & \multicolumn{2}{l}{$ a^*_{k+1} = e_{k+1}\T( B - U_{kk-1} M_{kk-1}^{-1} V_{k-1k-1}\T )   $} \\
    & \multicolumn{2}{l}{$ \beta_k, w_{k+1} = \textnormal{house}(a^*_{k+1})  $} \\
    & \multicolumn{2}{l}{$ W_k = \begin{bmatrix} W_{k-1} & w_{k+1} \end{bmatrix}  $} \\
\multicolumn{3}{l}{$\textnormal{end } $} \\    
\end{tabular} \\

To be clear, $U_{kk-1}$ denotes the matrix $U_{00}$
after $k$ left reflectors and $k-1$ right ones. It is defined as 
$ U_{kk-1} = \begin{bmatrix} \pl{b} & Y_{k} & BW_{k-1} \end{bmatrix} $. Similarly for
$ V_{kk-1} = \begin{bmatrix} \pl{c} & B\T Y_{k} & W_{k-1} \end{bmatrix} $ and
$ M_{kk-1} $.

The main computations in Algorithm 1 are the products $ a_k = A_{k-1k-1} e_k $
and $ a^*_{k+1} = e\T_{k+1} A_{kk-1} $. Not all elements of these vectors have to be computed, as we only need the last $m-k+1$ elements of $ a_k $ and the last $n-k$ of $ a^*_{k+1} $ for the reflectors. To analyze the dominant number of multiplications in the algorithm we assume that $ n \le m $ and that 
$ Y_k, W_k, R_k, T_k, Y_k \T b $ and $ c\T W_k $ are stored (see \cref{fig:mid}). We explicitly estimate the cost of $ a_k(k:m) = A_{k-1k-1}(k:m,:)e_k $, which in turn is similar to the cost of $ a^*_{k+1}(k+1:n) $.
To estimate $ a_k(k:m) $, note that $ V_{k-1k-1}\T e_k $ is not expensive;
the only operation that needs any flops is $ Y_{k-1} \T B e_k $ at
$2(k-1)$ multiplications. To solve with the middle matrix, two triangular 
solves are needed at $ \frac{k(k+1)}{2} + \frac{(k-1)k}{2} + k-1 = \mathcal{O}(k^2) $ multiplications, and the product $ Y_{k-1}\T B W_{k-1} v  $ with
some $k$ vector $v$. Using the nonzeros in $ Y_{k-1}, B $ and $W_{k-1}$ we may form
the product $ Y_{k-1}\T B W_{k-1} v  $ in 
$ 2n(k-1) + 2n - (k-1)^2 = \mathcal{O}(2nk - k^2)  $ multiplications. Combining
the previous terms shows that solving with the middle matrix needs
$\mathcal{O}(2nk - k^2 + k^2 = 2nk)$ flops. Next, the product $ U_{k-1k-1}(k:m,:) u $ for some $ 2k - 1 $ vector $u$ is $ \mathcal{O}((m+n)k - k^2) $.
Note that $ B e_k $ does not need any multiplications. Finally, combining the dominant operations, we form the product $ a_k(k:m) =  A_{k-1k-1}(k:m,:)e_k $ in 
\begin{equation*}
    \mathcal{O} \left((m+3n)k - 2k^2 \right)
\end{equation*}
flops.
As the complexity of $ e_{k+1} \T A_{kk-1}(:,k+1:n) = a^*_{k+1}(k+1:n) $ is 
similar to that of $ a_k(k:m)$, we conclude that one iteration of Algorithm 1 uses about 
\begin{equation*}
    \mathcal{O}\left(2((m+3n)k - 2k^2) \right)
\end{equation*}
flops. Running the complete bidiagonalization, i.e., for $k=1,\ldots,n-1$, therefore has complexity
\begin{equation}
    \label{eq:housecplx}
    \mathcal{O} \left((m+3n)n^2 - \frac{4}{3}n^3 \right) = \mathcal{O} \left(mn^2 + \frac{5}{3}n^3 \right).
\end{equation}
By storing the dense $k \times k$ array $ Y\T_k B W_k $ we may reduce the complexity to $ \mathcal{O} \left(mn^2 + \frac{2}{3}n^3 \right) $. Comparing \eqref{eq:housecplx} to the dense Householder method from \cref{sec:house}
\eqref{eq:cplxhouse}, we note that the compact representation has a flop count that is higher by approximately $2n^3$, although it does maintain the sparsity of the bidiagonal. This approach is therefore most advantageous if a partial factorization is useful or if $ n \ll m $ and memory is a concern.

\subsection{Givens Low Rank}
\label{sec:giv}
Even though Algorithm 1 can be good for certain situations, it essentially recomputes a bidiagonal
factorization (albeit preserving sparsity). 
To exploit
the special structure of 
$ B + \pl{b} {\pl{c}}\T $ we introduce a method based on a sequence of Givens rotations (see e.g., Golub and van Loan \cite{golub2013matrix}). We represent the orthogonal rotation that acts on either rows $i$ and $j$ 
(when applied from the left) or on columns $i$ and $j$ (when applied from the right) as a sparse rank-2 update
to the identity:
\begin{equation}
\label{eq:giv}
G^{\textnormal{d}}_{ij} = I + g^{\textnormal{d}}_1 e\T_i + g^{\textnormal{d}}_2 e\T_j, \quad \quad \textnormal{direction}=\textnormal{d} \in \{\textnormal{left},\textnormal{right} \},
\end{equation}

\begin{figure}   
\centering
\includegraphics[scale=0.64,trim={4.5cm 5cm 3.5cm 2cm},clip]{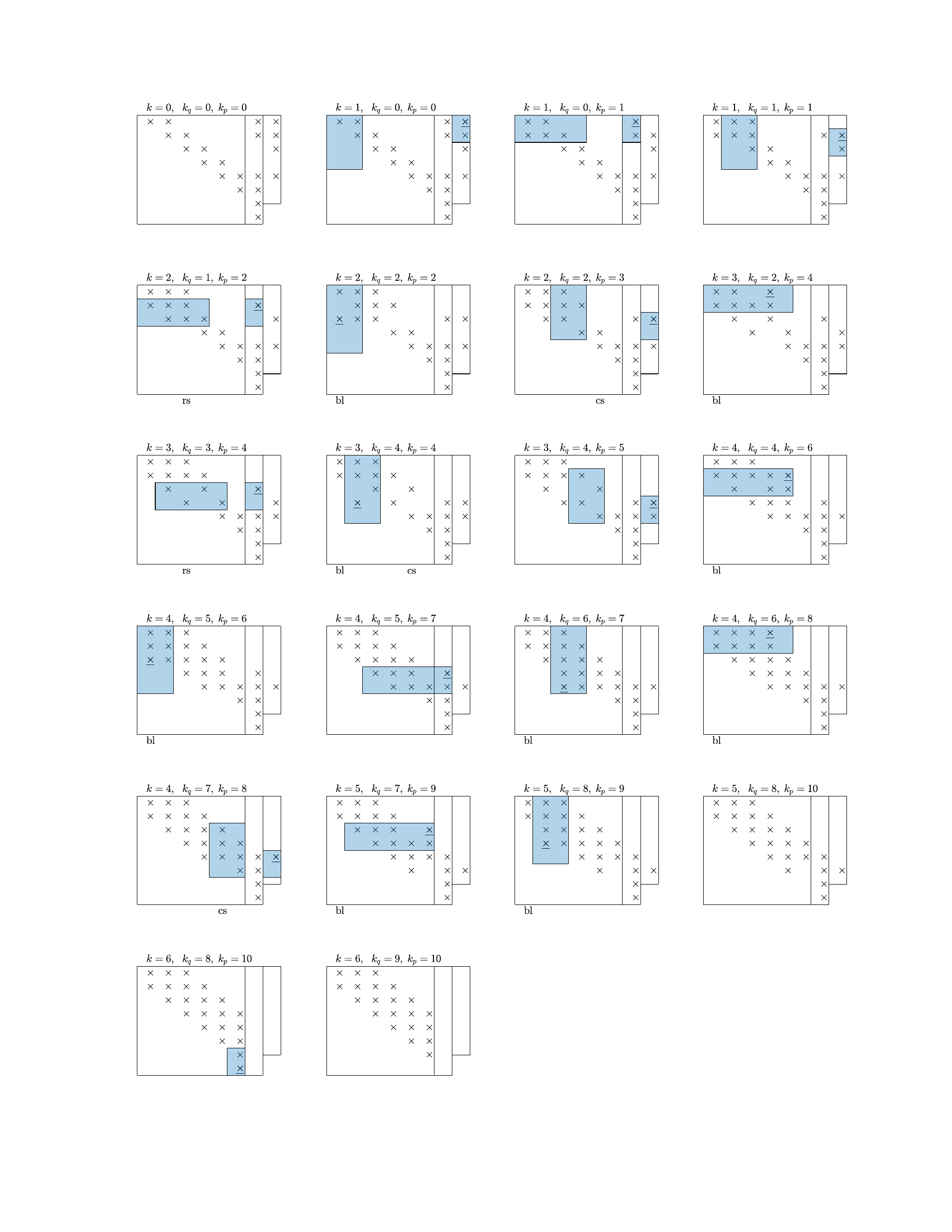}%
\caption{Phase 1 in a Givens bidiagonal updating algorithm. Shaded elements are accessed at a given iteration.
An underlined element is about to be eliminated. Labels rs, cs, and bl indicate if an iteration is a row or a column
permutation or a bulge chase. $ k_q $, $k_p$ count the number of left or right rotations} 
\label{fig:pdupp1}
\end{figure}

\vspace{-20pt}

\begin{equation}
    \label{eq:giv1}
    \left( g^{\textnormal{d}}_1 \right)_l = 
    \begin{cases}
        c^{\textnormal{d}} & \textnormal{ if } l = i \\
        s^{\textnormal{d}} & \textnormal{ if } l = j \\
        0                   & \textnormal{ otherwise }
    \end{cases} \quad
    \left( g^{\textnormal{d}}_2 \right)_l = 
    \begin{cases}
        -s^{\textnormal{d}} & \textnormal{ if } l = i \\
        c^{\textnormal{d}} & \textnormal{ if } l = j \\
        0                   & \textnormal{ otherwise }
    \end{cases}
    \quad (c^{\textnormal{d}})^2 + (s^{\textnormal{d}})^2 = 1.
\end{equation}
The actual values of $ c^{\textnormal{d}} = c $ and $ s^{\textnormal{d}} = s $ are typically determined so that
$ \big[ \begin{smallmatrix} c & -s \\ s & c \end{smallmatrix} \big] \big[ \begin{smallmatrix} \gamma_i \\ \gamma_j \end{smallmatrix} \big] = \big[ \begin{smallmatrix} 0 \\ \alpha \end{smallmatrix} \big]  $, where 
$ \gamma_i, \gamma_j $ are the $ i^{\textnormal{th}} $ and $ j^{\textnormal{th}} $ elements in a row or column
that is to be eliminated. 
Note that when $ \pl{b} = \gamma_n e_n + \gamma_{n-1} e_{n-1} $ and $ \pl{c} = \bar{\gamma}_n e_n $ 
for some scalars $ \gamma_n, \gamma_{n-1}, \bar{\gamma}_n $ then $ B + \pl{b} {\pl{c}}\T $ is in bidiagonal form and no further reduction
is needed. Subsequently our strategy alternately eliminates the elements in $ \pl{b} $ and $\pl{c}$ while maintaining
a certain sparsity in $B$. To be specific, in phase 1 we eliminate all necessary elements in $\pl{b}$ and $\pl{c}$ and
transform $B$ into a banded matrix $ B_{pq} $ with one nonzero subdiagonal $(p=1)$ and two nonzero superdiagonals $(q=2)$.
This matrix has 4 nonzero diagonals and therefore 2 more than the desired updated bidiagonal $\pl{B}$ with $p=0$, $q=1$.
In a second phase (phase 2) we reduce the quaddiagonal to bidiagonal. In contrast to \eqref{eq:bfill}
we illustrate that the sparsity pattern after applying one rotation is only minimally disturbed:
\begin{equation}
    \label{eq:gfill}
   \extsup{G}{l}{12} ( B + \pl{b} {\pl{c}} \T ) =
   \left[
    \begin{array}{c c c c }
    \times    & \times              &             &  \\
    \times    & \ddots            & \ddots            &  \\
                &             &           \times  & \times \\
                &             &         &               \times
    \end{array}
    \right] +
    \left[
    \begin{array}{c}
        0 \\
        \times \\
        \vdots \\
        \times
    \end{array}
    \right]
    \left[
        \begin{array}{c c c c}
            \times & \times & \hdots & \times 
        \end{array}
    \right].
\end{equation}
There are practical considerations for an efficient implementation of the algorithm though.
For instance, when two consecutive elements in $ \pl{b} $ or $ \pl{c} $ are $ \gamma_k \ne 0 $ and
$ \gamma_{k+1} = 0 $ then we permute these elements and the corresponding rows or columns of $B$.
Further, for $m > n$ the final $\pl{b}$ is typically not $ \pl{b} = \gamma_n e_n + \gamma_{n-1} e_{n-1} $ but
$ \pl{b} = \gamma_n e_n + \gamma_{n-1} e_{n-1} + [0; \pl{b}(n+1:m) ]\T $, i.e., all elements up to the $(n-1)^{\textnormal{th}}$ are zero,
with $m-n$ trailing nonzero elements. To eliminate this ``spike'' we apply a sequence of reflectors in
reverse order $ G_{n+1n} \cdots G_{m m-1} \pl{b}(n+1:m)  $, which importantly does not affect any elements of $B$
because this part is all zero. On the other hand, applying a rotation typically introduces an unwanted
nonzero even in a quaddiagonal matrix. We therefore have to ``chase'' these bulge elements off the matrix
by new sequences of rotations. To give some intuition on phase 1 of the algorithm, we show the nonzero patterns
of $B$, $\pl{b}$ and $\pl{c}$ after $k$ iterations in \cref{fig:pdupp1}. The labels $\textnormal{bl}$, $\textnormal{cs}$ and
$\textnormal{rs}$ below each pattern indicate whether this iteration is a bulge chase, a column permutation, or a row permutation.
An underline below an element means this element will be zeroed next. The shaded part shows which elements
are being accessed by the algorithm. 


Importantly, each rotation only needs $10=2\times5 $ multiplications (the maximum number of elements in the 
shaded areas in \cref{fig:pdupp1}). That is, each rotation is an $ \mathcal{O}(1) $ operation.
This is the basis for our argument of the complexity of this method. Specifically, there are up to $n-1$ elements
in each of $ \pl{b} $ and $\pl{c}$ that 
cause a bulge in $B$. A bulge in turn requires at most $n$ rotations
to be chased-off. Hence there are approximately $2n^2$ multiplications to eliminate
nonzeros in $\pl{b}$ and $\pl{c}$ and bring $B$ into banded from with bandwidth $p=1$ and $q=2$.
Thus, phase 1 requires $ \mathcal{O}(2n^2) $ flops. Phase 2 reduces the quaddiagonal to bidiagonal, which requires
chasing an additional $2n-3$ bulges.
The total complexity of $ \mathcal{O}(4n^2) $ 
is an order of magnitude less complex than regular bidiagonalization methods. We call this algorithm
\begin{center}
    Algorithm 2: BGU (Bidiagonal Givens Update)
\end{center}

\subsection{Implementation}
\label{sec:impl}
Both BHU (Algorithm 1) and BGU (Algorithm 2) are implemented in Fortran 90, with C interfaces to 
Matlab and Python \cite{bidiagSoftware}. The methods are compiled with gfortran from gcc 10.5.0 on a Linux computer with Ubuntu 22.04
and on a MacBook with gcc 14.2.0 and Sequoia 15.3 OS. We use the compiler optimization flag \texttt{-ofast}. To call the functions from Matlab we use the MEX interoperation layer, and to use them within Python we use ctypes from the standard library. The memory layout for BHU is similar to
\cref{fig:mid}, while for BGU we store each rotation as a vector of 4 elements, i.e., $c,s$ and indices $i,j$.

\section{Applications} 
\label{sec:numex}
Because of its quadratic complexity, algorithm BGU is suited for tracking streaming subspaces.
A particular application is recommendation systems in which a set of user ratings for specific items 
evolve over time. Another is network link prediction, in which a low-rank graph representation
is used to infer the likelihood of two users being connected in a (social) network. Other machine-learning tasks such as object tracking \cite{lim2004incremental} or fine-tuning \cite{han2023svdiff} as well as optimization \cite{brustTOMSSR1} may also be formulated in a streaming context. To update a low-dimensional subspace,
suppose that at a given time 
the factors $ Q^{(1)} \in \mathbb{R}^{m \times r} $, $ P^{(1)} \in \mathbb{R}^{n \times r} $ and $ B \in \mathbb{R}^{r \times r} $ are known and updated with $ b \in \mathbb{R}^m $ and $c \in \mathbb{R}^n $:
\begin{equation}
\label{eq:lrup}
\pl{A} = Q^{(1)} B {P^{(1)}}\T + b c\T.
\end{equation}
Define projections of $b$ and $c$ onto the orthogonal complement of 
$ Q^{(1)} $ and $ P^{(1)} $ as
\begin{equation}
\label{eq:defup}
    \pl{b} := {Q^{(1)}}\T b, \quad b^{\perp} := b - Q^{(1)} \pl{b}, \quad
    \pl{c} := {P^{(1)}}\T c, \quad c^{\perp} := c - P^{(1)} \pl{c},  
\end{equation}
\vspace{-1cm}
\begin{equation*}
    \delta :=  \| b^{\perp}  \|_2, \quad \gamma := \| c^{\perp} \|_2.
\end{equation*}
Further, assume for the moment that $ \delta > 0 $ and $\gamma > 0 $, and note that 
the QR factorization of $ \bmat{Q^{(1)} & b } $ is given by 
\begin{equation}
    \label{eq:qrlr}
    \bmat{Q^{(1)} & b } = \bmat{Q^{(1)} & b^{\perp} / \delta }  \bmat{ I_r & \pl{b} \\ & \delta } := \bar{Q}^{(1)} \bar{R}^{(1)}.    
\end{equation}
The QR factorization of $ \bmat{P^{(1)} & c } $ is analogous with interchanges $ Q^{(1)} \leftrightarrow P^{(1)}$, $ \pl{b} \leftrightarrow \pl{c}$,
$ b^{\perp} \leftrightarrow c^{\perp} $ and $ \delta \leftrightarrow \gamma $ in \cref{eq:qrlr}.
The low-rank update \cref{eq:lrup} 
becomes
\begin{equation}
    \label{eq:equivlrup}
    \pl{A} = \bmat{Q^{(1)} & b^{\perp} / \delta } 
    \left( \bmat{ B & \\ & 0 } + \bmat{\pl{b} \\ \delta } \bmat{ {\pl{c}}\T & \gamma  } \right) \bmat{ P^{(1)} & c^{\perp} / \gamma }\T.
\end{equation}
We may now apply Algorithm 2 BGU to the $(r+1)\times (r+1)$ middle matrix of 
\cref{eq:equivlrup} for an updated $(r+1)\times (r+1)$ bidiagonal factorization 
$ Q^{(2)} \pl{\bar{B}} {P^{(2)}}\T $. In the notation from
\cref{eq:qrlr}, an exact rank $r+1$ factorization of $\pl{A}$ is given by
\begin{equation}
    \label{eq:bdlrup}
    \pl{A} = \bar{Q}^{(1)} Q^{(2)} \pl{\bar{B}} (\bar{P}^{(1)} P^{(2)})\T.
\end{equation}
The rank-$r$ factorization is now obtained by removing the row and column of
$ \pl{\bar{B}} $ that correspond to the lowest column-norm index in $ \pl{\bar{B}} $. Without loss of generality, assume that the last column norm is smallest,
so that the rank-$r$ matrix is defined by $ \pl{\bar{B}}(1:r,1:r) $.
Note that for large problems, the most expensive computations are the products
$ \bar{Q}^{(1)}(:,:) Q^{(2)}(:,1:r)  $ and $ \bar{P}^{(1)}(:,:) P^{(2)}(:,1:r) $. If, on the other hand, $ \delta = 0 $ and/or $ \gamma=0 $ then 
the columns of $ \bar{Q}^{(1)} $ and/or $ \bar{P}^{(1)} $ are not augmented by 
$ \pl{b} $ and/or $ \pl{c} $, and the factorization of \cref{eq:equivlrup} is
automatically rank-$r$.
We apply our algorithms
to large link prediction and recommendation problems and benchmark them against state-of-the-art methods.

\subsection{Link Prediction}
\label{sec:linkpred}
To predict a connection between users $i$ and $j$ in a network represented by an adjacency matrix, the estimate $ e\T_i A e_j $
may be used. At time $h>0$, the update
\begin{equation}
    \label{eq:lrpred}
        A_h = A_{h-1} + b_h c_h\T
\end{equation}
often corresponds to adding a new user to the network. This situation can be modeled through 
$ A_0 = 0_{n \times n} $ and $ b_h \in \{0,1\}^n $, where $ b_h(i)=1 $ indicates that users $i$ and $h$ are 
connected and not otherwise, and $ c_h = e_h $. Deng et al.~\cite{dengIncSVD} recently developed updating methods
for a truncated SVD for large and sparse systems. In their numerical experiments the authors compare with methods
by Zha and Simon \cite{zha1999updating}, Vecharynski and Saad \cite{vecharynski2014fast} and Yamazaki et al.~\cite{yamazaki2015randomized},
and conclude that their version of RPI (randomized power iteration) is most effective. We use the public repository
of this study's implementations for comparison. We also use the same datasets and a similar test setup \cite{dengIncSVD}. 

The Flickr data \cite{Zafarani+Liu:2009} 
is a friendship network among bloggers on Flickr, an image and video hosting website.
This network has 80,513 users and 5,899,882 friendship pairs, which corresponds to 11,799,764 edges in the graph. Slashdot \cite{leskovec2009community}
is a technology-related news website known for its specific user community. The network contains 82,168 users 
and 948,464 edges. 

\begin{figure}
    \centering
    \begin{tikzpicture}
        \matrix[matrix of nodes, row sep=0pt, column sep=0.0cm] {
            \node {
                \begin{tikzpicture}
                    \begin{axis}[
                        width=7.5cm,
                        height=6cm,
                        xlabel={$r$},
                        ylabel={Seconds},
                        title style={font=\small},
                        ylabel style={font=\small},
                        tick label style={font=\scriptsize},
                        grid=none,
                        ymax=2000,
                        ymin=-399,
                        xtick={2048,2560,3072,3584,4096,4608,5120},
                        ytick={0,500,1000,1500},
                        xticklabel style={rotate=0},
                        title={Flickr (80,513 nodes, 11,799,764 edges)},
                        thick,
                        clip=false
                    ]
                    \addplot+[
                        mark=*,
                        color=blue,
                        nodes near coords,
                        point meta=y,
                        every node near coord/.append style={
                            font=\scriptsize,
                            /pgf/number format/fixed,
                            /pgf/number format/precision=1,
                            yshift=-15pt
                        }
                    ] coordinates {
                        (2048, 167.085247)
                        (2560, 244.557526)
                        (3072, 337.697366)
                        (3584, 432.572201)
                        (4096, 555.201138)
                        (4608, 692.367875)
                        (5120, 829.957809)
                    };
                    \addplot+[
                        mark=triangle*,
                        color=red,
                        nodes near coords,
                        point meta=y,
                        every node near coord/.append style={
                            font=\scriptsize,
                            /pgf/number format/fixed,
                            /pgf/number format/precision=1,
                            yshift=15pt
                        }
                    ] coordinates {
                        (2048, 167.971788)
                        (2560, 252.103524)
                        (3072, 377.640275)
                        (3584, 539.749382)
                        (4096, 796.714689)
                        (4608, 1072.857309)
                        (5120, 1414.555905)
                    };
                    \addplot+[
                        only marks,
                        mark=none,
                        nodes near coords,
                        point meta=explicit symbolic,
                        every node near coord/.append style={
                            font=\tiny,                            
                            yshift=-23pt,
                            text=blue
                        }
                    ] coordinates {
                        (2048, 167.085247) [1e-11] 
                        (2560, 244.557526) [1e-11]
                        (3072, 337.697366) [5e-12]
                        (3584, 432.572201) [9e-12]
                        (4096, 555.201138) [9e-12]
                        (4608, 692.367875) [1e-11]
                        (5120, 829.957809) [7e-12]
                    };
                    \addplot+[
                        only marks,
                        mark=none,
                        nodes near coords,
                        point meta=explicit symbolic,
                        every node near coord/.append style={
                            font=\tiny,
                            yshift=7pt,
                            text=red
                        }
                    ] coordinates {
                        (2048, 167.971788) [1e-12]
                        (2560, 252.103524) [3e-12]
                        (3072, 377.640275) [3e-12]
                        (3584, 539.749382) [5e-12]
                        (4096, 796.714689) [5e-12]
                        (4608, 1072.857309) [4e-12]
                        (5120, 1414.555905) [5e-12]
                    };
                    \end{axis}
                \end{tikzpicture}
            }; 
            &
            \node {
                \begin{tikzpicture}
                    \begin{axis}[
                        width=7.5cm,
                        height=6cm,
                        xlabel={$r$},
                        title style={font=\small},
                        ylabel style={font=\small},
                        tick label style={font=\scriptsize},
                        grid=none,
                        ymax=2000,
                        ymin=-399,
                        xtick={2048,2560,3072,3584,4096,4608,5120},
                        ytick={0,500,1000,1500},
                        xticklabel style={rotate=0},
                        title={Slashdot (82,168 nodes, 948,464 edges)},
                        thick,
                        clip=false
                    ]
                    \addplot+[
                        mark=*,
                        color=blue,
                        nodes near coords,
                        point meta=y,
                        every node near coord/.append style={
                            font=\scriptsize,
                            /pgf/number format/fixed,
                            /pgf/number format/precision=1,
                            yshift=-15pt
                        }
                    ] coordinates {
                      (2048, 169.734196)
                      (2560, 251.899696)
                      (3072, 345.775095)
                      (3584, 445.855686)
                      (4096, 569.759965)
                      (4608, 715.083963)
                      (5120, 961.168424)
                    };
                    \addplot+[
                        mark=triangle*,
                        color=red,
                        nodes near coords,
                        point meta=y,
                        every node near coord/.append style={
                            font=\scriptsize,
                            /pgf/number format/fixed,
                            /pgf/number format/precision=1,
                            yshift=15pt
                        }
                    ] coordinates {
                      (2048, 157.969878)
                      (2560, 237.93629)
                      (3072, 347.157336)
                      (3584, 482.89094)
                      (4096, 740.989401)
                      (4608, 978.398822)
                      (5120, 1319.079309)
                    };
                    \addplot+[
                        only marks,
                        mark=none,
                        nodes near coords,
                        point meta=explicit symbolic,
                        every node near coord/.append style={
                            font=\tiny,
                            yshift=-23pt,
                            text=blue
                        }
                    ] coordinates {
                        (2048, 169.730000) [4e-12]
                        (2560, 251.900000) [1e-11]
                        (3072, 345.780000) [6e-12]
                        (3584, 445.860000) [7e-12]
                        (4096, 569.760000) [1e-11]
                        (4608, 715.080000) [1e-11]
                        (5120, 961.170000) [7e-12]
                    };
                    \addplot+[
                        only marks,
                        mark=none,
                        nodes near coords,
                        point meta=explicit symbolic,
                        every node near coord/.append style={
                            font=\tiny,
                            yshift=7pt,
                            text=red
                        }
                    ] coordinates {
                        (2048, 157.970000) [6e-13]
                        (2560, 237.940000) [1e-12]
                        (3072, 347.160000) [1e-12]
                        (3584, 482.890000) [8e-13]
                        (4096, 740.990000) [1e-12]
                        (4608, 978.400000) [1e-12]
                        (5120, 1319.100000) [1e-12]
                    };
                    \end{axis}
                \end{tikzpicture}
            }; \\
        };
    \end{tikzpicture}

    \vspace{0.5em}
    \begin{tikzpicture}
        \begin{axis}[
            hide axis,
            xmin=0, xmax=1,
            ymin=0, ymax=1,
            legend columns=2,
            legend style={
                at={(0.5,1.05)},
                anchor=south,
                draw=none,
                /tikz/every even column/.append style={column sep=1cm}
            }
        ]
        \addlegendimage{mark=*,blue}
        \addlegendentry{BGU (our)}
        \addlegendimage{mark=triangle*,red}
        \addlegendentry{Deng et al.~\cite{dengIncSVD} }
        \end{axis}
    \end{tikzpicture}

    \caption{Comparison of algorithm BGU and randomized power iteration (RPI)~\cite{dengIncSVD} on two link prediction datasets. Each data point shows the total time of all updates and the cumulative absolute error below. The errors are similar and close to machine precision, while BGU scales well for increasing rank.}
    \label{fig:linkpred}
\end{figure}
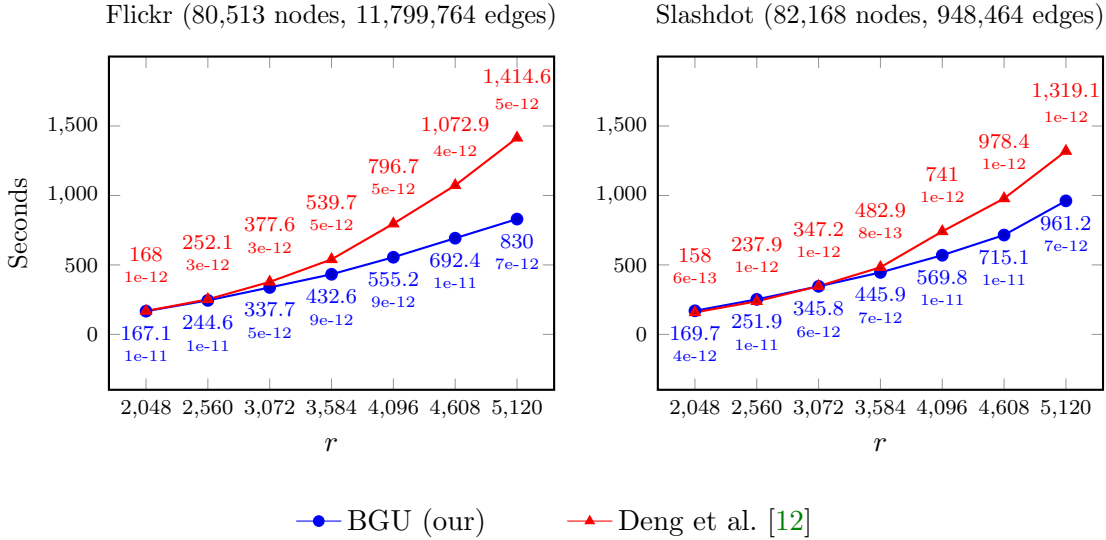

The RPI method updates a rank-$r$ SVD with factors $ U^{(1)}, \Sigma, V^{(1)} $, and has two additional
integer parameters $t$ and $l$. Parameter $l$ is called the ``spatial dimension of the approximation'' and is set to $l=10$ as in the article.
Parameter $t$ corresponds to ``subiterations'' and is set to $t=3$ as in \cite{dengIncSVD}. To avoid runtime errors
we set the sparse flag to false. This experiment demonstrates the scaling of BGU and RPI (the best of multiple previous methods)
for increasing sizes of $r$, namely $r = 256\cdot \{7,9,11,13,15,17,19\}$. The stream starts from zero: 
\begin{equation}
    \label{eq:stream}
        A_0 = 0_{n \times n}, \quad U^{(1)}_0 = Q^{(1)}_0 = V^{(1)}_0 = P^{(1)}_0 = I_n(:,1:r), \quad \Sigma_0 = B_0 = 0_{r\times r},
\end{equation}
and is run for 50 updates: $ h = 1,2,\ldots 50 $. We compute the error in fitting the true network as 
\begin{equation}
\label{eq:res}
    \sum_{h=1}^{50} \left| \| A_h \|_F - \| \Sigma_h \|_F \right| \quad \textnormal{ or } \quad
                         \sum_{h=1}^{50} \left| \| A_h \|_F - \| B_h \|_F \right|.
\end{equation}
The outcomes are in \cref{fig:linkpred}.

\subsection{Recommendation Systems}
\label{sec:recsys}
\begin{figure} 
    \centering    
    \begin{tikzpicture}
        \node[align=center] at (0,0) {MovieLens 32M  \\  (87,585 movies, 200,948 users, 32,000,204 ratings)};
    \end{tikzpicture}
    \begin{tikzpicture}
        \matrix[matrix of nodes, row sep=0pt, column sep=0cm] {
            \node {
                \begin{tikzpicture}
                    \begin{axis}[ 
                        width=6.9cm,    
                        height=6cm,     
                        xlabel={Update},
                        ylabel={Residuals},
                        title style={font=\footnotesize},
                        ylabel style={font=\footnotesize},
                        xlabel style={font=\footnotesize},
                        tick label style={font=\scriptsize},
                        ymode=log,                        
                        xticklabel style={rotate=0},
                        thick,
                        clip=false
                    ]
                    \addplot[
                        no markers,
                        thick,
                        color=blue,
                        ] table [x index=0, y index=15] {data/bdup_r2000.txt};
                    \addplot[
                        no markers,
                        thick,
                        color=red,
                        ] table [x index=0, y index=12] {data/brand_r2000.txt};
                    \end{axis}
                \end{tikzpicture}
            }; 
            &
            \node {
                \begin{tikzpicture}
                \begin{axis}[
                        axis y line*=right,
                        width=6.9cm, 
                        height=6cm,  
                        xlabel={Update},
                        ylabel={Seconds/Update},
                        ylabel style={font=\footnotesize},
                        xlabel style={font=\footnotesize},
                        tick label style={font=\scriptsize},
                        thick,
                        clip=false
                    ]
                    \addplot[
                        only marks,
                        mark=*,
                        mark size=0.25pt,
                        color=blue,
                        opacity=0.05,
                    ] table [x index=0, y index=3] {data/bdup_r2000.txt}; 
                    \addplot[
                        only marks,
                        mark=triangle,
                        mark size=0.25pt,
                        color=red,
                        opacity=0.05,
                    ] table [x index=0, y index=3] {data/brand_r2000.txt}; 
                \end{axis}
                \begin{axis}[
                        axis y line*=left,
                        axis x line=none,
                        width=6.9cm, 
                        height=6cm,  
                        ylabel={Seconds (total)},
                        ylabel shift = -5pt,
                        ylabel style={font=\footnotesize},
                        tick label style={font=\scriptsize},                        
                        thick,
                        clip=false
                    ]
                    \addplot[
                        no markers,
                        thick,
                        color=blue    
                    ] table [x index=0, y index=16] {data/bdup_r2000.txt};
                    \addplot[
                        no markers,
                        thick,
                        color=red,
                    ] table [x index=0, y index=13] {data/brand_r2000.txt};
                \end{axis}
                \end{tikzpicture}
            }; \\
        };
    \end{tikzpicture}

    \begin{tikzpicture}
        \begin{axis}[            
            hide axis,
            xmin=0, xmax=1,
            ymin=0, ymax=1,
            legend columns=2,
            legend style={
                at={(0.5,1.05)},
                anchor=south,
                draw=none,
                /tikz/every even column/.append style={column sep=1cm}
            }
        ]
        \addlegendimage{blue} 
        \addlegendentry{BGU (our)}
        \addlegendimage{red} 
        \addlegendentry{Brand \cite{brand2006fast}}
        \end{axis}
    \end{tikzpicture}
    \caption{Algorithm BGU and the incremental SVD \cite{brand2006fast} for updating a rank-$r=2000$
    approximation. The residuals remain small in both methods. Each BGU update costs about 0.25 seconds,
    while the incremental SVD takes about 0.5 secs (blue and red dots in the right axis). 
    }
    \label{fig:recsys}
\end{figure}
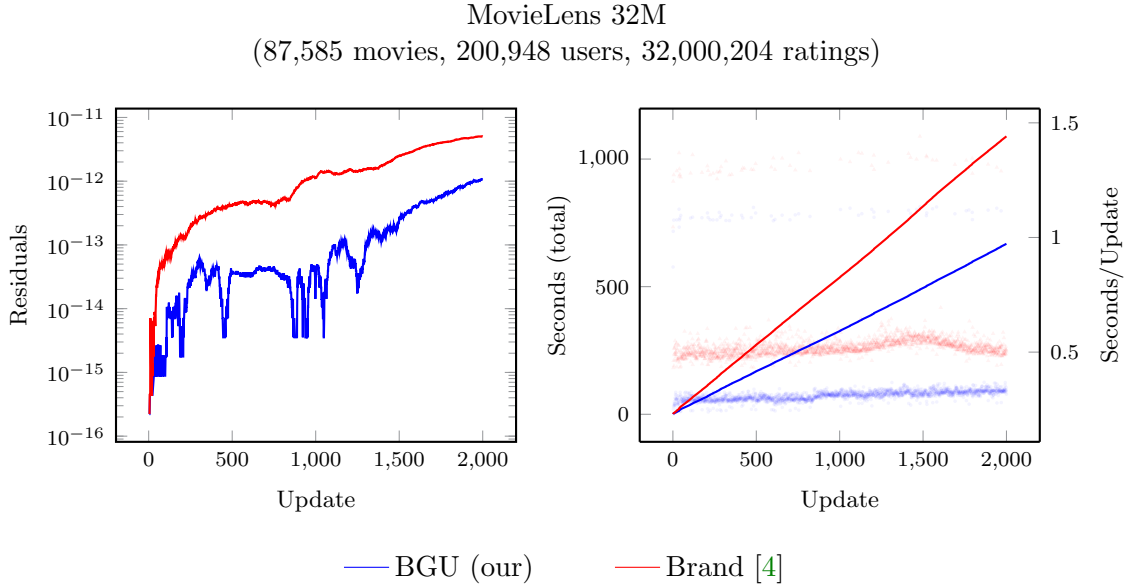

Consider a database of movies, say, each mapped to an index in the set
$ \mathcal{I} = \{1,\ldots,m \} $.
At time $ h = 0,1,\ldots $ a specific movie selection is $ i_{h} \in \mathcal{I} $.
Further, consider a set of users mapped to indices $ \mathcal{J} = \{1, \ldots, n \} $. 
At a given time, user $ j_h \in \mathcal{J} $ may choose to rate movie $i_h$ with rating $ \theta_{i_h j_h} \in \{1,2,3,4,5 \} $,
where $1$ is lowest and $5$ is highest score. If all ratings up to time $h$ are arranged into 
matrix $A_h \in \mathbb{R}^{m \times n} $, the stream $ A_h, A_{h+1}, \ldots $ represents the evolution
of movie preferences. We assume that one time point is associated with one particular preference, i.e., $ \theta_{i_h j_h} = \theta_h $.
If two preferences were submitted at exactly the same time, the tie is split randomly so
that one is first and one second. This means that the streaming system can be modeled as
\begin{equation}
    \label{eq:streamsys}
    A_{h+1} = A_h + \theta_h e_{i_h} e_{j_h} \T, \quad A_0 = 0.
\end{equation}
The MovieLens 32M dataset \cite{harper2015movielens} contains
32,000,204 ratings from 200,948 users on 87,585 movies. Each user has rated at least 20 movies. The data was collected through the MovieLens website 
between January 09, 1995 and October 12, 2023. For comparison, we use the incremental SVD due to Brand \cite{brand2006fast}. 
This algorithm was shown to be faster than iterative methods, 
such as svds 
based on the ARPACK software \cite{lehoucq1998arpack}. Brand's
algorithm is similar to the way we apply BGU for subspace updating. However, instead of maintaining a bidiagonal factorization
it uses a truncated SVD with factors $ U^{(1)}_h, \Sigma_h, V^{(1)}_h  $, and when finding the next SVD of the 
corresponding middle matrix in \cref{eq:equivlrup} (i.e., after swapping $ U^{(1)}_h \leftrightarrow Q^{(1)}_h $, $ V^{(1)}_h \leftrightarrow P^{(1)}_h $
and $ \Sigma_h \leftrightarrow B_h $) it does not use an updating strategy but a complete rank-$(r+1)$ SVD. 
It is known that products such as $ \bar{Q}^{(1)} Q^{(2)} $ or $ \bar{U}^{(1)} U^{(2)} $ can deteriorate in orthogonality,
and that it may be necessary to reorthogonalize these matrices occasionally. We use a heuristic to check for loss of orthogonality
periodically, by monitoring its change. We found in our experiments with thousands of iterations that reorthogonalizing
only a few times may be effective. As in the previous experiment from \Cref{sec:linkpred} \cref{eq:stream} we initialize the stream from zero:
\begin{equation*}
        A_0 = 0_{m \times n}, \quad U^{(1)}_0 = Q^{(1)}_0 = I_m(:,1:r), \quad V^{(1)}_0 = P^{(1)}_0 = I_n(:,1:r), \quad \Sigma_0 = B_0 = 0_{r\times r},
\end{equation*}
and perform $h+1=1,2,\ldots,2000$ updates. We define the residual at each update as 
\begin{equation*}
    \left| \| A_h \|_F - \| \Sigma_h \|_F \right| \quad \textnormal{ or } \quad
                        \left| \| A_h \|_F - \| B_h \|_F \right|.
\end{equation*}
\Cref{fig:recsys} shows that algorithm BGU accurately fits the data and is faster than 
the incremental SVD.



\subsection{Benchmarking}
\label{sec:benchmark}
\begin{figure}
\centering
\includegraphics[trim={0.5cm 0.5cm 0.5cm 0.5cm},clip]{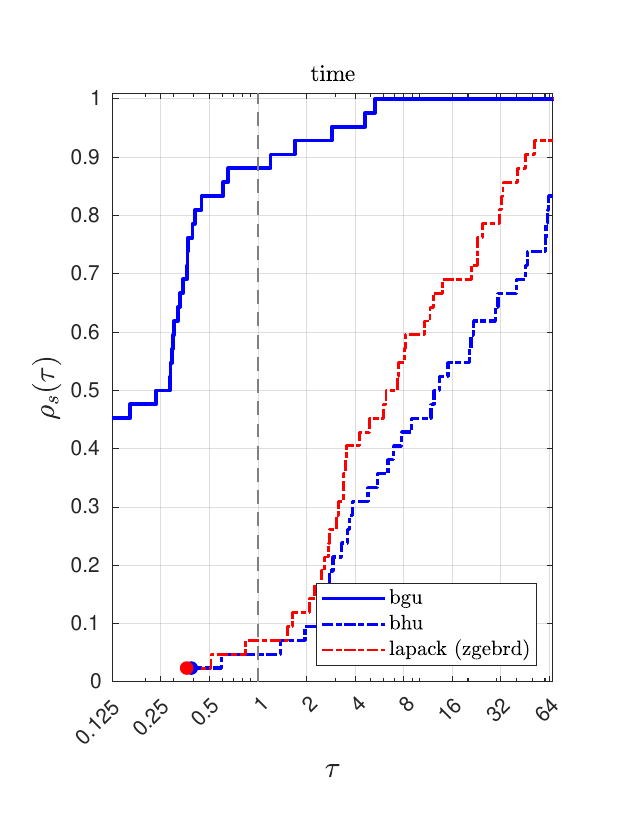}
\caption{Algorithms BGU, BHU and lapack's zgebrd on a set of 43 SuiteSparse problems. The data from \Cref{tab:exp1} is
displayed using a performance profile~\cite{dolan2002benchmarking} of the computational time 
($ \rho_s(\tau) $ is the fraction of problems for which solver $s$ was within a factor $\tau$ of the fastest solver). BGU has quadratic complexity and is significantly
faster then zgebrd.}
\label{fig:lapackcomp}
\end{figure}

In a third experiment we compare our proposed algorithms to the LAPACK algorithm zgebrd, which is 
considered a very effective implementation. Algorithm zgebrd is written in Fortran 77,
with interfaces to Matlab, Python and many other languages. We use a set of 43 SuiteSparse \cite{SuiteSparseMatrix} matrices
with dimensions  $ 500 \le m \le 15000 $ and $n$ less than 40\% of the number of rows: $ n < \lfloor 0.4 m \rfloor $.
After a rank-1 correction we compute the error $ |\|\pl{A}\|_F - \|\pl{B}\|_F| $. We only compute a thin bidiagonalization with 
zgebrd and let the algorithm determine its own ``optimal'' block size. The results of applying BHU, zgebrd and BGU are
shown in \cref{tab:exp1} and summarized in \Cref{fig:lapackcomp}.  For all methods, the computed errors are similar, and except for 5 relatively small problems, BGU is the fastest.
For certain larger problems, the differences become very apparent. For instance, on problem \texttt{ch7-6-b3} (a combinatorial matrix)
zgebrd takes about 76 seconds, while BGU requires about 1.5 seconds. For problem \texttt{Fran6}, the times are
about 22 seconds for zgebrd and 0.67 seconds for BGU.

\subsection{Future Directions}
\label{sec:future}
We have shown that the proposed algorithms are effective for updating the bidiagonal factorization. As this 
is the precursor to virtually all dense SVD algorithms, future work may extend this to an efficient
and accurate method for updating the SVD. For symmetric systems, the proposed updating
algorithms may be adapted to generate a tridiagonal, which is the main step to compute eigenvalues.
The subspace tracking algorithm may be sped-up by avoiding rectangular-square matrix products
when in \cref{eq:defup} $\delta = 0 $ or not (and the same for $\gamma$), and by refining the 
reorthogonalization strategy.
Further,
testing the new methods on fine-tuning and optimization applications appears to be a promising direction. 

\section{Conclusions}
\label{sec:conclusions}
We have proposed new methods to update a bidiagonalization following low-rank changes. 
The methods are well suited to streaming applications, where a large number of updates must be
processed quickly. Our Householder updating method is effective at maintaining sparsity and has
a lower memory footprint than typical Householder methods, at the expense of higher flops. However,
our Givens updating method is one order lower in complexity than standard algorithms, scales well, and can
be very fast and accurate for processing high volumes of low-rank data changes.

\appendix

\begin{table}[t]  
 \captionsetup{width=1\linewidth}
\caption{Comparison of BHU, BGU and lapack's zgebrd algorithms for updating a bidiagonal with rank-one change on 42 matrices from SuiteSparse \cite{SuiteSparseMatrix}. 
In column 3, ``Dty"  is the density of a particular matrix $ A $ calculated as $ \text{Dty} = \frac{\text{nnz}(A)}{m \cdot n} $.
Bold entries mark the
fastest times, while second fastest are italicized. The table's results are summarized in \cref{fig:lapackcomp}.
}
\label{tab:exp1}

\setlength{\tabcolsep}{2pt} 

\small
\hbox to 1.00\textwidth{\hss 
\begin{tabular}{|l | c | c | c c | c c | c c |} 
		\hline
		\multirow{2}{*}{Problem} & \multirow{2}{*}{$m / n$} & \multirow{2}{*}{\text{Dty}} & 
		\multicolumn{2}{c|}{BHU } &  		
		\multicolumn{2}{c|}{lapack (zgebrd)} &  	
        \multicolumn{2}{c|}{BGU } \\
		\cline{4-9}
		& & 
		&  Sec & Res
		&  Sec & Res
		&  Sec & Res \\
		\hline
$\texttt{ash608}$ & 608/188 & 0.01 &0.011&2e-16&\emph{0.0061}&0&\textbf{0.004}&0\\  
$\texttt{ash958}$ & 958/292 & 0.007 &0.034&2e-16&\emph{0.015}&0&\textbf{0.0063}&2e-16\\
$\texttt{illc1033}$ & 1033/320 & 0.01 &0.045&2e-16&\emph{0.02}&0&\textbf{0.0065}&6e-16\\ 
$\texttt{illc1850}$ & 1850/712 & 0.007 &0.44&0&\emph{0.13}&0&\textbf{0.037}&0\\ 
$\texttt{well1033}$ & 1033/320 & 0.01 &0.045&0&\emph{0.019}&0&\textbf{0.007}&1e-15\\ 
$\texttt{well1850}$ & 1850/712 & 0.007 &0.43&1e-16&\emph{0.12}&0&\textbf{0.035}&0\\ 
$\texttt{photogrammetry}$ & 1388/390 & 0.02 &0.087&1e-11&\emph{0.035}&1e-12&\textbf{0.011}&7e-11\\ 
$\texttt{photogrammetry2}$ & 4472/936 & 0.009 &2&1e-11&\emph{0.56}&1e-12&\textbf{0.069}&3e-10\\ 
$\texttt{Franz1}$ & 2240/768 & 0.003 &0.59&0&\emph{0.19}&0&\textbf{0.045}&3e-15\\ 
$\texttt{Franz5}$ & 7382/2882 & 0.002 &58&0&\emph{20}&0&\textbf{0.63}&8e-15\\ 
$\texttt{Franz6}$ & 7576/3016 & 0.002 &66&0&\emph{22}&0&\textbf{0.67}&1e-14\\ 
$\texttt{Franz7}$ & 10164/1740 & 0.002 &21&0&\emph{11}&0&\textbf{0.35}&6e-16\\ 
$\texttt{GL7d11}$ & 1019/60 & 0.02 &\emph{0.0067}&5e-16&\textbf{0.0034}&0&0.0098&0\\ 
$\texttt{GL7d12}$ & 8899/1019 & 0.004 &5.3&0&\emph{2.8}&0&\textbf{0.11}&1e-16\\ 
$\texttt{ch5-5-b2}$ & 600/200 & 0.01 &0.013&2e-16&\emph{0.006}&0&\textbf{0.0037}&2e-16\\ 
$\texttt{ch6-6-b2}$ & 2400/450 & 0.007 &0.18&0&\emph{0.074}&0&\textbf{0.012}&2e-16\\ 
$\texttt{ch7-6-b1}$ & 630/42 & 0.05 &\emph{0.00059}&0&0.0016&0&\textbf{0.00022}&0\\ 
$\texttt{ch7-6-b2}$ & 4200/630 & 0.005 &0.64&1e-16&\emph{0.25}&0&\textbf{0.031}&5e-16\\ 
$\texttt{ch7-6-b3}$ & 12600/4200 & 0.001 &2e+02&1e-16&\emph{76}&1e-16&\textbf{1.5}&7e-15\\ 
$\texttt{ch7-7-b1}$ & 882/49 & 0.04 &\emph{0.029}&2e-16&\textbf{0.01}&2e-16&0.056&0\\ 
$\texttt{ch7-7-b2}$ & 7350/882 & 0.003 &3.1&4e-16&\emph{1.6}&0&\textbf{0.068}&1e-15\\ 
$\texttt{ch7-8-b1}$ & 1176/56 & 0.04 &\textbf{0.0019}&0&0.0039&0&\emph{0.0032}&1e-16\\ 
$\texttt{ch7-8-b2}$ & 11760/1176 & 0.003 &9.7&3e-16&\emph{6.2}&2e-16&\textbf{0.15}&3e-15\\ 
$\texttt{ch7-9-b1}$ & 1512/63 & 0.03 &0.0085&0&\textbf{0.0061}&0&\emph{0.0073}&5e-16\\ 
$\texttt{ch8-8-b1}$ & 1568/64 & 0.03 &\emph{0.0027}&0&0.0055&0&\textbf{0.00092}&3e-16\\ 
$\texttt{cis-n4c6-b2}$ & 1330/210 & 0.01 &0.021&0&\emph{0.012}&1e-16&\textbf{0.00098}&0\\ 
$\texttt{cis-n4c6-b3}$ & 5970/1330 & 0.003 &7.1&2e-16&\emph{2.9}&0&\textbf{0.036}&0\\ 
$\texttt{mk10-b1}$ & 630/45 & 0.04 &\textbf{0.00069}&2e-16&\emph{0.0018}&0&0.0032&2e-16\\ 
$\texttt{mk10-b2}$ & 3150/630 & 0.005 &0.48&1e-16&\emph{0.17}&0&\textbf{0.016}&1e-16\\ 
$\texttt{mk11-b1}$ & 990/55 & 0.04 &\emph{0.0013}&2e-16&0.0029&0&\textbf{0.00039}&0\\ 
$\texttt{mk11-b2}$ & 6930/990 & 0.003 &3.9&0&\emph{1.9}&0&\textbf{0.042}&4e-16\\ 
$\texttt{mk12-b1}$ & 1485/66 & 0.03 &\emph{0.0031}&1e-16&0.0059&1e-16&\textbf{0.0012}&3e-16\\ 
$\texttt{mk12-b2}$ & 13860/1485 & 0.002 &19&0&\emph{13}&0&\textbf{0.17}&1e-15\\ 
$\texttt{mk9-b2}$ & 1260/378 & 0.008 &0.083&1e-16&\emph{0.032}&0&\textbf{0.0093}&0\\ 
$\texttt{n2c6-b3}$ & 1365/455 & 0.009 &0.12&0&\emph{0.044}&0&\textbf{0.002}&0\\ 
$\texttt{n3c6-b3}$ & 1365/455 & 0.009 &0.12&2e-16&\emph{0.042}&0&\textbf{0.002}&0\\ 
$\texttt{n4c5-b3}$ & 1350/455 & 0.009 &0.12&0&\emph{0.043}&0&\textbf{0.0031}&2e-16\\ 
$\texttt{n4c6-b2}$ & 1330/210 & 0.01 &0.021&0&\emph{0.012}&0&\textbf{0.001}&0\\ 
$\texttt{n4c6-b3}$ & 5970/1330 & 0.003 &7&0&\emph{2.8}&0&\textbf{0.033}&0\\ 
$\texttt{rel6}$ & 2340/157 & 0.01 &0.021&4e-16&\emph{0.016}&0&\textbf{0.0057}&2e-16\\ 
$\texttt{relat6}$ & 2340/157 & 0.02 &0.02&2e-16&\emph{0.015}&0&\textbf{0.0043}&0\\ 
$\texttt{abtaha1}$ & 14596/209 & 0.02 &0.26&4e-16&\emph{0.15}&0&\textbf{0.068}&9e-16\\
 \hline
 \end{tabular}
 \hss}
\end{table} 

\clearpage

\section*{Acknowledgments}
We are very grateful for Tamara Kolda's valuable feedback on our original manuscript, which prompted comparing with two additional algorithms on realistic problems.

\bibliographystyle{siamplain}
\bibliography{refs}

\end{document}


\maketitle

\section{Bidiagonal Truncation Error}
\label{sec:AppBidiagtrErr}
We show how to establish the error bound in \cref{eq:errtrBD}. A related result can be found
in \cite{simon2000low}. First, without loss of generality, assume that $A$ is square and
consider the $n \times n$ orthogonal projection matrices $E_r$ and $\bar{E}_r$:
\begin{equation}
    \label{eq:Er}
    E_r := 
    \begin{bmatrix}
        I(:,1\!:\!r) & 0(:,r\!+\!1\!:\!n)
    \end{bmatrix},
    \quad 
    \bar{E}_r :=
    \begin{bmatrix}
         0(:,1\!:\!r) & I(:,r\!+\!1\!:\!n) 
    \end{bmatrix},
    \quad 
    E_r + \bar{E}_r = I.
\end{equation}
Note that $\extsup{A}{BD}{r}$ from \cref{eq:trBD} can be expressed in matrix 
notation using $E_r$.  In particular,
\begin{equation}
    \label{eq:trErrFact}
    \| A - \extsup{A}{BD}{r} \|^2_F = \| A - Q E_r B E\T_r P\T \|^2_F.
\end{equation}
By carefully multiplying $E_r B E_r\T$ we find that 
\begin{equation}
    \label{eq:product}
Q E_r B E_r\T P\T = Q (B  - \bar{E}_r B \bar{E}_r\T  - \beta_r e_r e_{r+1}\T )P\T.
\end{equation}
Using the bidiagonal decomposition of $A$ from \cref{eq:bidiag}, substituting 
\cref{eq:product} into \cref{eq:trErrFact} and using the property that the Frobenius
norm is invariant under orthognal transformations, we obtain
\begin{equation*}
    \| A - \extsup{A}{BD}{r} \|^2_F = \| \bar{E}_r B \bar{E}_r\T  + \beta_r e_r e_{r+1}\T \|^2_F.
\end{equation*}
There is a special form for the matrix on the right, namely
\begin{equation*}
    \bar{E}_r B \bar{E}_r\T  + \beta_r e_r e_{r+1}\T =
    \left[
    \begin{array}{c c c | c c c}
        0   &           & & & & \\
            & \ddots    &   &           & & \\
            &           & 0 & \beta_r       & & \\
\hline      &           &   & \alpha_{r+1}  & \ddots & \\
            &           &   &   & \ddots & \beta_{n-1} \\
            &           &   &   &   & \alpha_{n}
    \end{array}
    \right].
\end{equation*}
Therefore the norm can be evaluated exactly as
\begin{equation*}
    \| A - \extsup{A}{BD}{r} \|^2_F = \| \bar{E}_r B \bar{E}_r\T  + \beta_r e_r e_{r+1}\T \|^2_F = \sum_{i=r+1}^n \alpha^2_i + \beta^2_{i-1}.
\end{equation*}

\section{Between Truncation Error}
\label{sec:AppBetweenErr}
We use the 
definitions in \cref{eq:Er}.
As the SVD and the bidiagonal decomposition of $A$ both exist, there must be two orthogonal matrices
$ Q_1 $ and $P_1$ for which
\begin{equation*}
    Q_1\T B P_1 = \Sigma.
\end{equation*}
Therefore a second upper bound is 
\begin{equation*}
    \| \extsup{A}{SVD}{r} - \extsup{A}{BD}{r} \|^2_F \le \| \extsup{A}{SVD}{r} \|^2_F + \| \extsup{A}{BD}{r} \|^2_F =
    \| E_r Q_1\T  B  P_1 E_r\T \|^2_F + \| Q_1\T E_r B E_r\T P_1 \|^2_F.
\end{equation*}
Because 
$$
    \| E_r Q_1\T  B  P_1 E_r\T \|^2_F = \sum_{i=1}^r \sigma^2_i 
\qquad \mathrm{and} \qquad 
\| Q_1\T E_r B E_r\T P_1 \|^2_F =  \sum_{i=1}^r \alpha^2_i + \beta^2_{i-1},
$$
 we see that 
 \begin{equation*}
    \| \extsup{A}{SVD}{r} - \extsup{A}{BD}{r} \|^2_F \le \sum_{i=1}^r \sigma^2_i + \alpha^2_i + \beta^2_{i-1}.
\end{equation*}
In fact, $ \| \extsup{A}{SVD}{r} - \extsup{A}{BD}{r} \|^2_F $ can be evaluated exactly to provide the lower bound.
Without loss of generality, assume that $t=n$ with  
$ P_1 = \begin{bmatrix} {p^*_1}\T & \cdots & {p^*_r}\T & \cdots {p^*_t}\T \end{bmatrix}\T $ given row-wise. Then
\begin{equation}
    \label{eq:exactbound}
    \| \extsup{A}{SVD}{r} - \extsup{A}{BD}{r} \|^2_F = \| Q_1 \Sigma_r P_1\T - B_r \|^2_F =
    \sum_{i=1}^r \| Q_1 \Sigma_r {p^*_i}\T - b_i \|^2_2 + \sum_{i=r+1}^t \| Q_1 \Sigma_r {p^*_i}\T \|^2_2.
\end{equation}
The $i^{\textnormal{th}}$ term in the first sum on the right is 
\begin{equation}
    \label{eq:exactbound1}
    \| Q_1 \Sigma_r {p^*_i}\T - b_i \|^2_2 = \| \Sigma_r {p^*_i}\T \|^2_2 - 2 p^*_i (Q_1 \Sigma_r) \T b_i + \| b_i \|^2_2. 
\end{equation}
Because 
\begin{equation}
    \label{eq:exactbound2}
    \Sigma_r = Q_1 \T B P_1 - \sum_{j=r+1}^t \sigma_j e_j e_j\T, \quad (Q_1 e_j)\T b_i = \sigma_j p^*_i e_j,
\end{equation}
we evaluate the term in the middle of \eqref{eq:exactbound1} as
\begin{equation}
    \label{eq:exactbound3}
    p^*_i (Q_1 \Sigma_r) \T b_i = \| b_i \|^2_2 - \sum_{j=r+1}^t \sigma^2_j (p^*_i e_j)^2.
\end{equation}
Combining all terms in \eqref{eq:exactbound} using \eqref{eq:exactbound1} and \eqref{eq:exactbound3} we find that
\begin{equation}
    \label{eq:exactbound4}
    \| \extsup{A}{SVD}{r} - \extsup{A}{BD}{r} \|^2_F =
    \sum_{i=1}^r \sigma^2_i - (\alpha^2_i + \beta^2_{i-1}) + 2  \sum_{j=r+1}^t \sum_{i=1}^r \sigma^2_j  (p^*_i e_j)^2.
\end{equation}
Observe that the last term in \eqref{eq:exactbound4} is always nonnegative; thus a lower bound is given by
\begin{equation}
    \label{eq:exactlower}
    \sum_{i=1}^r \sigma^2_i - (\alpha^2_i + \beta^2_{i-1}) \le \| \extsup{A}{SVD}{r} - \extsup{A}{BD}{r} \|^2_F.
\end{equation}
An upper bound can also be found from the exact expression \eqref{eq:exactbound4} and the fact that $ \sum_{i=1}^r (p^*_i e_j)^2 \le 1  $ because it is a
partial squared length of an orthogonal vector. However, the resulting bound is equal to the one given in the text below \eqref{eq:diffErr}.

\section{Representation of the bounds in \cref{eq:diffErr} }


\begin{figure}[h]   
\centering
\includegraphics[scale=0.95,trim={0.8cm 0.8cm 0.8cm 1cm},clip]{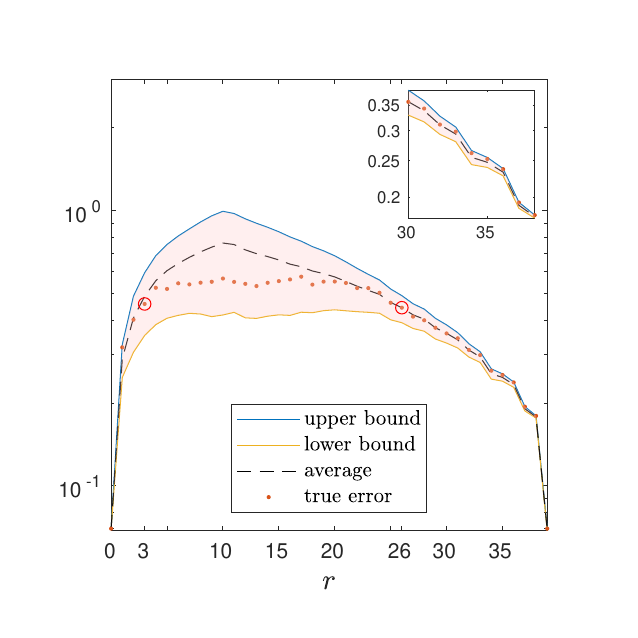}%
\caption{Bounds for the normalized difference between a rank-$r$ SVD or BD approximation on a $60 \times 40$ matrix. The true error is nonlinear, and can be similar even for $r$ values that are far apart (the two circles at $r=3$ and $r=26$). Overall, the bounds are tight, especially at both ends of the rank spectrum.}
\label{fig:bounds_supp}
\end{figure}   

\Cref{fig:bounds_supp} illustrates that the bounds between the truncated BD and SVD are tight, especially 
close to the endpoints of the rank spectrum.




 







\bibliographystyle{siamplain}
\bibliography{refs}